\documentclass[12pt]{amsart}
\usepackage{amsmath,amscd,amssymb,amssymb}

\begin{document}

\numberwithin{equation}{section}

\theoremstyle{plain}
\newtheorem{theorem}{Theorem}[section]
\newtheorem{prop}[theorem]{Proposition}
\newtheorem{lemma}[theorem]{Lemma}
\newtheorem{cor}[theorem]{Corollary}
\newtheorem{obs}[theorem]{Observation}
\newtheorem{conj}[theorem]{Conjecture}
\newtheorem{question}[theorem]{Question}
\newtheorem{example}[theorem]{Example}
\newtheorem{claim}[theorem]{Claim}

\theoremstyle{definition}
\newtheorem{defn}[theorem]{Definition}
\newtheorem{remark}[theorem]{Remark}
\newtheorem{defnprop}[theorem]{Proposition-Definition}

\newcommand{\Z}{{\mathbb Z}}
\newcommand{\N}{{\mathbb N}}
\newcommand{\OO}{{\mathfrak{o}}}

\newcommand{\HH}{{\mathcal H}}
\newcommand{\II}{{\mathcal I}}
\newcommand{\UU}{{\mathcal U}}
\newcommand{\FF}{{\mathcal F}}
\newcommand{\RR}{{\mathcal R}}
\newcommand{\EE}{{\mathcal E}}

\newcommand{\Oa}{{\OO}_1}
\newcommand{\Ob}{{\OO}_2}
\newcommand{\Oas}{{\OO}_1^{\times}}
\newcommand{\Obs}{{\OO}_2^{\times}}

\newcommand{\ib}{{\mathbf i}}
\newcommand{\rb}{{\mathbf r}}

\newcommand{\sig}{{\sigma}}

\newcommand{\pl}{{\delta}}
\newcommand{\pll}{{\pi^\sig}}

\newcommand{\DD}{{\mathcal D}}
\newcommand{\CC}{{\mathcal C}}
\newcommand{\PP}{{\mathcal I}}

\newcommand{\ha}{{\hat{0}}}
\newcommand{\hb}{{\hat{1}}}

\newcommand{\jc}{{l}}

\newcommand{\T}{{\DD}}

\newcommand{\x}{M_{\mu}}
\newcommand{\xx}{M_{\mu'}}

\newcommand{\y}{M_{\la}}
\newcommand{\z}{M_{\la/\mu}}

\newcommand{\m}{\mathfrak{p}}
\let\myacute=\'
\let\mygrave=\`
\let\myddot=\"

\newcommand{\hlf}{\genfrac{}{}{0.1pt}{1}{1}{2}}

\newcommand{\la}{{\lambda}}
\newcommand{\gb}{{\textbf{g}}}

\newcommand{\f}[1]{\flat{#1}}

\newcommand{\G}{{G}}
\newcommand{\K}{{K}}
\newcommand{\J}{{J}}

\newcommand{\E}{{E}}
\newcommand{\Pp}{{P}}
\newcommand{\Pt}{\mathcal{P}}

\newcommand{\GL}{\text{GL}}
\newcommand{\SL}{\text{SL}}
\newcommand{\B}{\text{B}}
\newcommand{\U}{\text{U}}
\newcommand{\Conj}{\text{Conj}}

\newcommand{\End}{\text{End}}
\newcommand{\Stab}{\text{Stab}}
\newcommand{\Hom}{\text{Hom}}
\newcommand{\Aut}{\text{Aut}}
\newcommand{\Ind}{\text{Ind}}
\newcommand{\Ker}{\text{Ker}}
\newcommand{\Res}{\text{Res}}
\newcommand{\Cusp}{\text{Cusp}}
\newcommand{\Inf}{\text{Inf}}
\newcommand{\Inv}{\text{Inv}}

\newcommand{\F}{\mathbb{F}}
\newcommand{\Fb}{\mathbf{F}}

\newcommand{\lx}{\ell_1}
\newcommand{\ly}{\ell_2}
\newcommand{\kk}{\ell}
\newcommand{\ls}{\ell}

\newcommand{\lxx}{\lx'}
\newcommand{\lyy}{i}

\newcommand{\Ok}{{\OO}_\kk}
\newcommand{\Oks}{{\OO}_\kk^{\times}}
\newcommand{\Oj}{{\OO}_j}
\newcommand{\Ojs}{{\OO}_j^{\times}}

\newcommand{\ldiff}{d \la}

\newcommand{\lu}{l}
\newcommand{\ld}{l'}

\newcommand{\ue}[2]{#1 \hookrightarrow #2}
\newcommand{\uq}[2]{#1 \twoheadrightarrow #2}

\newcommand{\ie}[2]{\mathbf{i}_{#1, \hookrightarrow}^{#2}}
\newcommand{\iq}[2]{\mathbf{i}_{#1,\twoheadrightarrow}^{#2}}
\newcommand{\re}[2]{\mathbf{r}_{#1,\hookrightarrow}^{#2}}
\newcommand{\rqq}[2]{\mathbf{r}_{#1,\twoheadrightarrow}^{#2}}

\newcommand{\is}[2]{\mathbf{i}_{#1}^{#2}}
\newcommand{\rs}[2]{\mathbf{r}_{#1}^{#2}}

\newcommand{\M}{M_2(A)}
\newcommand{\val}{v}
\newcommand{\ord}{o}

\newcommand{\Fq}{\mathbb{F}_q}

\newcommand{\Fqx}{\mathbb{F}_q^{\times}}
\newcommand{\C}{\mathbb{C}}
\newcommand{\inv}{^{-1}}
\newcommand{\bsl}{\backslash}
\newcommand{\half}{\frac{1}{2}}
\newcommand{\mat}[4]{\left(\begin{array}{cc}#1 & #2 \\ #3 & #4\end{array}\right)}
\newcommand{\kappabar}{\overline{\kappa}}
\newcommand{\lt}[1]{\langle #1 \rangle}

\newcommand{\gfrac}[2]{\genfrac{}{}{0.1pt}{1}{#1}{#2}}


\title{Representations of automorphism groups \linebreak of finite {\large $\OO$}-modules of rank
two}

\author{Uri Onn }

\subjclass[2000]{20G05, 22E50}

\thanks{Supported by the Israel Science Foundation, ISF
grant no. 555104, by the Edmund Landau Minerva Center for Research
in Mathematical Analysis and Related Areas, sponsored by the Minerva
Foundation (Germany) and by Marie Curie research training network
LIEGRITS (MRTN-CT 2003-505078).}

\maketitle

\begin{abstract}
Let $\OO$ be a complete discrete valuation domain with finite
residue field. In this paper we describe the irreducible
representations of the groups $\Aut(M)$ for any finite $\OO$-module
$M$ of rank two. The main emphasis is on the interaction between the
different groups and their representations. An induction scheme is
developed in order to study the whole family of these groups
coherently. The results obtained depend on the ring $\OO$ in a very
weak manner, mainly through the degree of the residue field. In
particular, a uniform description of the irreducible representations
of $\GL_2(\OO/\m^\kk)$ is obtained, where $\m$ is the maximal ideal
of $\OO$.

\end{abstract}


\section{Introduction}

\subsection{Overview} Let $F$ be a non-Archimedean local field with ring of integers
$\OO=\OO_{F}$, maximal ideal $\m=\m_F$, a fixed uniformizer
$\pi=\pi_F$ and residue field of cardinality $q$. Let $\tilde{\OO}
\supset \OO$ be the ring of integers of a degree $n$ non-ramified
extension $\tilde{F} \supset F$ and let $\tilde{\OO}_{\kk} \supset
\OO_{\kk}$ be the reductions modulo $\m^{\kk}$ ($\kk \in \N$). Let
$\Lambda_n$ denote the set of partitions of length $n$ and
$\Lambda=\cup \Lambda_n$. For any $\la=(\ell_1,\ldots,\ell_n) \in
\Lambda$ let $M_{\lambda}$ denote the $\OO$-module
$\oplus_{i}\OO_{\ell_i}$ of type $\la$. Endow $\Lambda$ with a
partial order where $\mu \le \la$ if and only if a module of type
$\mu$ can be embedded in a module of type $\la$. Let
\[
\G_{\lambda}=\Aut_{\OO}(M_{\lambda}).
\]

In this paper the complex irreducible representations of the groups
$G_{\la}~(\la \in \Lambda_2)$ are classified and discussed in
detail. The Harish-Chandra philosophy of cusp forms \cite{HC} is
adopted and implemented in order to interconnect the representation
theories of these groups. Representations of $\G_{\la}$ are built
from representations of $\G_{\mu}~(\mu < \la)$ using various
induction functors. Representations which cannot be obtained from
groups of lower type are called cuspidal and are constructed as
well.

This approach was carried out by Green \cite{Green} who classified
the representations of the groups $\G_{1^n}=\Aut_{\Fq}(\Fq^n)$ $(n
\in \N)$. The cuspidal representations of $\G_{1^n}$, i.e., those
which are not contained in parabolically induced representations
from $\G_{1^m}$ $(m < n)$, are parameterized by an object of
arithmetic flavor: orbits of primitive characters of
$\F_{q^n}^{\times}$ under the Galois action. On the other hand, the
induction part of the theory is of combinatorial flavor, is
independent of the field, and essentially coincides with the
representation theory of the symmetric group, see \cite{Zelevinsky}.
Such a picture is the model we seek in general. In this paper we
achieve it for automorphism groups of rank two $\OO$-modules.

The special types of rectangular shape $\kk^n=(\kk,\ldots,\kk)$ draw
much of attention since the representations of the groups
$\G_{\kk^n}$ $(\kk \in \N)$ exhaust the complex continuous
representations of $\GL_n(\OO_{F})$, the maximal compact subgroup of
$\GL_n(F)$. One of the main goals of this paper is to provide
evidence that the groups $\G_{\la}$ ($\la < \kk^n$) play
indispensable role in the classification of irreducible
representations of $\G_{\kk^n}$.

\subsection{The main results} We begin with a short description of the
subgroups of $\G_{\la}$ which play the role of maximal parabolics
for the groups $G_{1^n}~(n >1)$. The first family of these subgroups
is a straightforward generalization. Given a direct sum
decomposition $M_{\la}=M_{\mu} \oplus M_{\la/\mu}$ ($\la/\mu$
denotes the type of the quotient $M_{\la}/M_{\mu}$), let
$\Pp_{\mu,\la}$ be the stabilizer of $M_{\mu}$ in $M_{\la}$. Then
$\Pp_{\mu,\la}$ surjects on $\G_{\mu} \times \G_{\la/\mu}$, and one
can construct representations of $G_{\la}$ by pulling back
representations from the product $\G_{\mu} \times \G_{\la/\mu}$ to
$\Pp_{\mu,\la}$ and inducing them to $\G_{\la}$. Since the term \lq
parabolic\rq~ is quite loaded we shall call the groups
$\Pp_{\mu,\la}$ {\em (maximal) geometric stabilizers} in $\G_{\la}$.

The second family has no analogue in the case of the groups
$G_{1^n}$. Given types $\mu \le \la$, call $\mu$ {\em symmetric} in
$\la$ if the embedding of $M_{\mu}$ in $M_{\la}$ is unique up to
automorphism. In such case let $\la / \mu$ denote the well defined
type of the quotient $M_{\la}/M_{\mu}$. As will be explained later
on, it follows that there is a unique embedding (up to automorphism)
of a module $M_{\la/\mu}$ of type $\la/\mu$ such that
$M_{\la}/M_{\la/\mu}$ is of type $\mu$. Let $\mu$ be a symmetric type
in $\la$. Assume also that as partitions $\mu$ and $\la$ are of same
length, or equivalently, that the corresponding modules have the
same rank. This assumption guarantees that there is no overlap with
the previous family. Set
\begin{align*}
&\Pp_{\mu \hookrightarrow \la}=\Stab_{\G_{\la}}(M_{\mu} \subset
M_{\la}), & &\\
&\Pp_{\la  \twoheadrightarrow \mu}=\Stab_{\G_{\la}}(M_{\la/\mu}
\subset M_{\la}).& &
\end{align*}

The fact that $\mu$ is symmetric in $\la$ translates into
canonical epimorphisms from $\Pp_{\mu \hookrightarrow \la}$ and
$\Pp_{\la \twoheadrightarrow \mu}$ to $\G_{\mu}$. The first is
defined by restriction of elements of $\Pp_{\mu \hookrightarrow
\la}$ to $M_{\mu}$, and the second by identifying elements of
$\Pp_{\la \twoheadrightarrow \mu}$ as automorphisms of
$M_{\la}/M_{\la/\mu}$. Again, we construct representations of
$\G_{\la}$ by pulling back representations from $\G_{\mu}$ to
either $\Pp_{\mu \hookrightarrow \la}$ or $\Pp_{\la
\twoheadrightarrow \mu}$, and then inducing them to $\G_{\la}$. We
shall call the groups $\Pp_{\la \twoheadrightarrow \mu}$ and
$\Pp_{\mu \hookrightarrow \la}$ {\em infinitesimal stabilizers} in
$\G_{\la}$. The types which eventually contribute to the
infinitesimal induction are types $\mu$ which are symmetric in
$\la$ and of same length and same height (the latter stands for the largest part of the partition). We denote the set of these
types by $\PP_{\la}$.

By a {\em twist} of a representation we mean its tensor product with
a one dimensional representation. For a partition $\la \in
\Lambda_n$ let $\f{\la}=(\ell_1-1,\ldots,\ell_n-1)$. There is a
natural reduction map
\begin{equation*}
 \G_{\la}=\Aut_{\OO}(M_{\la})
\twoheadrightarrow \G_{\f{\la}}=\Aut_{\OO}(\pi M_{\la})
\end{equation*}
defined by restriction of automorphisms of $M_{\la}$ to $\pi
M_{\la}$ (this generalizes the usual reduction modulo $\m^{\kk-1}$
in the case of $\la=\kk^n$). Call a representation $\rho \in
\hat{\G}_{\la}$ {\em primitive} if none of its twists is a pullback
from $\hat{\G}_{\f{\la}}$. By using an inductive approach one can
focus only on primitive representations as the non-primitive
representations are twists of representations of groups of lower
type.

\begin{defn}
A primitive irreducible representation of $\G_{\la}$ is {\em
cuspidal} if none of its twists is contained in some geometrically
or infinitesimally induced representation. The set of cuspidal
representations in $\hat{\G}_{\la}$ is denoted by $\hat{C}_{\la}$.
\end{defn}

\begin{theorem}\label{maintheorem} For any $\la \in \Lambda_2$ the primitive irreducible representations of $\G_{\la}$
fall precisely into one of the following classes
\begin{enumerate}
\item Contained in a geometrically induced representation.

\item Infinitesimally induced from a unique cuspidal representation
$\rho \in \hat{C}_{\mu}$ ($\mu \in \PP_{\la}$) up to twist.

\item Cuspidal.
\end{enumerate}
\end{theorem}

The cuspidal representations $\hat{C}_{\la}$ ($\la \in \Lambda_2$)
are explicitly constructed in section
\ref{sec:constructioncuspidal}. In \cite{AOPS} it is further shown
that there is a canonical bijection between cuspidal representations
of $\G_{\kk^n}$ and Galois orbits of strongly primitive characters of
$\tilde{\OO}_{\kk}^{\times}$ whenever $n$ is prime.

\medskip
Let $f_m=f_m^{\la,\OO}$ denote the number of irreducible
representations of $G_{\la}$ ($\la \in \Lambda$) of dimension $m$.

\begin{conj}\label{conj2} For any $\la \in \Lambda$,

(Weak version) the isomorphism type of the group algebra
$\C\G_{\la}$ depends only on $\la$ and $q=|\OO/\m|$.

(Strong version) $f_m=f_m^{\la} \in \mathbb{Q}[q]$.
\end{conj}

\begin{theorem}\label{poly}Let $R_{\la,\OO}(\T)=\sum_m  f_m \T^m \in \Z[\T]$. Then
\begin{align*}
&(1)~R_{\la}(\T)=q^{\ell-2}(q-1)^2\T+q^{\ell-2}(q^2-1)\T^{q-1}+q^{\ell-2}(q-1)^3\T^q,&
 &\la=(\ell,1), \ell > 1.& \\
&(2)~R_{\la}(\T)=q
R_{\f{\la}}(\T)+q^{\ell_1+\ell_2-3}(q^2-1)\T^{q^{\ell_2-1}(q-1)}+q^{\ell_1+\ell_2-3}(q-1)^3\T^{q^{\ell_2}},&
~~~&\la=(\ell_1,\ell_2), \ell_1>\ell_2>1.& \\
&(3)~R_{\ell^2}(\T)=q
R_{(\ell-1)^2}(\T)+\gfrac{1}{2}(q-1)(q^2-1)q^{2\ell-3}\T^{q^{\ell-1}(q-1)}+\\
&\qquad \qquad \qquad \qquad \qquad
q^{2\ell-2}(q-1)\T^{q^{\ell-2}(q^2-1)}+ \gfrac{1}{2} q^{2
\ell-3}(q-1)^3\T^{q^{\ell-1}(q+1)},&  &\ell > 1.&
\end{align*}
In particular, the strong version of Conjecture \ref{conj2} holds
for all $\la \in \Lambda_2$.
\end{theorem}

The dimensions of the irreducible representations seem to be
polynomials in $q$ as well. Moreover, we conjecture that for $\la =
\kk^n$ the following holds.

\begin{conj}\label{conjdim} The dimensions of primitive irreducible
representations of $\G_{\kk^n}$ are polynomials in $\mathbb{Z}[q]$ of degree
$d$ with $ d \le \left({n \atop
2}\right)\kk$.
\end{conj}

This holds for example when $\la=1^n$ by Green's work (in this case
we treat all the representations as primitive) or when $\la=\kk^2$
by Theorem \ref{poly}.

\subsection{History of the problem} A special case of the problem under consideration, namely when the
group is specialized to the general linear group
$\GL_2(\OO_{\kk})=\Aut_{\OO}(\OO_{\kk}^2)$ or to its closely related
subgroup, the special linear group $\SL_2(\OO_{\kk})$, has been
intensively studied by several authors. First results concerning the
representation theory of the group $\SL_2(\Z/p^{\kk}\Z)$ were
obtained by Kloosterman in his Annals papers
\cite{Kloosterman2,Kloosterman1}. Kloosterman constructed some of
the representations of $\SL_2(\Z/p^{\kk}\Z)$ using the Weil
representation (see also \cite{Springer} for a description of
Kloosterman's work). Tanaka \cite{Tanaka} generalized Kloosterman's
construction and obtained all the irreducible representations of
$\SL_2(\Z/p^{\kk}\Z)$. All these results assume that $p$ is odd. In
a series of papers by Nobs and Nobs-Wolfart
\cite{Nobs1,Nobs2,Nobs3,NobsWolfart1,NobsWolfart2}, these
constructions are generalized in several directions: some of the
constructions are applied for $p=2$, some of them are valid when
$\Z_p$ is replaced by any ring of integers of a local field, and
some are applied to $\GL_2$.  Other publications which do not use the Weil representation
include: Kutzko \cite{Kutzko} and  Nagornyi \cite{Nagornyi1} for
$\SL_2(\Z/p^{\kk}\Z)$, Silberger \cite{Silbereger2} for
$\text{PGL}_2(\OO_{\kk})$ and Shalika \cite{Shalika} for
$\SL_2(\OO_{\kk})$ when $\OO$ has characteristic zero. In a more recent paper \cite{JZ}, Jaikin-Zapirain computes the representation
zeta function for $\SL_2(\OO)$ in the case of odd characteristic. Nevertheless, only in a very recent preprint \cite{stasinski-2008} by Stasinski all the
irreducible representations of $\GL_2(\OO_{\kk})$ are constructed for general
$\OO$. For higher
dimensions there are some partial results of a more general nature
by Hill \cite{HG1,HG2,HG3,HG4}, Lusztig \cite{Lusztig1} and
Aubert-Onn-Prasad-Stasinski \cite{AOPS} however, the general case still
remains out of reach.

\subsection{Contents of the paper} The paper is organized as
follows. In section \ref{sec:operations} the various induction and
restriction functors are discussed in complete generality and the
concept of cuspidality is introduced. In section \ref{sec:rank2} we
specialize to the rank two case and explain the structure of the
various automorphism groups involved and their stabilizing
subgroups. In section \ref{sec:Gk1} we settle the case of the groups
$\G_{(\kk,1)}$ ($\kk>1$) due to their special structure and role.
The classification of their irreducible representations is simpler
comparing to the other automorphism groups. In section
\ref{sec:inductionscheme} we explain the inductive scheme in detail.
In section \ref{sec:constructioncuspidal} we construct the cuspidal
representations $\hat{C}_{\la} \subset \hat{\G}_{\la}$ for $\la$
with $\lx > \ly
> 1$ and briefly review the construction of $\hat{C}_{\kk^2}$ which
appeared in \cite{Nobs1} (but we shall follow the route taken in
\cite{AOPS}). Section \ref{sec:inductions} is devoted to geometric
and infinitesimal inductions. In section \ref{sec:unitarydual} we
collect all the results to obtain the classification of irreducible
representations of $\G_{\la}$ ($\la \in \Lambda_2$), in particular,
Theorems  \ref{maintheorem} and \ref{poly} are proved. In Section
\ref{sec:conjclasses} we classify the conjugacy classes in
$\G_{\la}$ ($\la \in \Lambda_2$).

\subsection{Acknowledgment} My deepest gratitude to Amritanshu Prasad and Anne-Marie Aubert for a most enjoyable
collaboration which was the seed for this paper. Many thanks to
Yakov Varshavsky and Alex Lubotzky for hosting the research, and to
David Kazhdan, Bernhard Keller and Dipendra Prasad for fruitful discussions. It is also a pleasure to thank Alexander Stasinski for his careful reading
of the paper and constructive comments. Finally, I thank the referee for valuable remarks.


\section{Operations on representations and cuspidality}\label{sec:operations} For a type
$\la \in \Lambda$ let $\CC_{\la}$ be the category of finite
dimensional complex representations of $\G_{\la}$. We shall now
define a family of functors between the categories $\{\CC_{\la} \}$
that generalize the usual functors of parabolic induction and their
adjoints which were defined by Green in \cite{Green} for the family
$\{\CC_{1^n}\}$ (see also \cite{Zelevinsky}).

\subsection{Geometric induction and
restriction}\label{geometricinduction}

The first family is the natural generalization of \cite{Green} which
comes from direct sum decompositions. Let $\y=\x \oplus \xx$ be a
decomposition of a finite $\OO$-module $\y$ into a direct sum of
$\OO$-modules of types $\mu$ and $\mu'$. Let $\Pp_{\x,\y} \subset
\Aut_{\OO}(\y)$ be the stabilizing subgroup of $\x$ in $\y$. Then we
have a surjection
\begin{equation}\label{geominduction}
\Aut_{\OO}(\y) \supset \Pp_{\x,\y} \overset{\iota}{\to}
\Aut_{\OO}(\x) \times \Aut_{\OO}(\y/\x),
\end{equation}
where the left part of $\iota$ is restriction and the right part
comes from the action of $\Pp_{\x,\y}$ on $\y/\x$. The functors of
geometric induction and geometric restriction are defined by
\begin{align*}
& & \is{\mu, \mu'}{\la}&: \CC_{\mu} \times \CC_{\mu'} \to \CC_{\la}
   & &(\xi,\xi') \mapsto \Ind_{\Pp_{\x,\y}}^{\Aut_{\OO}(\y)}(\iota^{*}(\xi \boxtimes \xi')
), &\qquad (\xi,\xi') \in \CC_{\mu} \times \CC_{\mu'}& & & \\
& & \rs{\mu, \mu'}{\la}&: \CC_{\la} \to \CC_{\mu} \times \CC_{\mu'}
& &\eta \mapsto \iota_{*}\Res_{\Pp_{\x,\y}}^{\Aut_{\OO}(\y)}(\eta),
&\qquad \eta \in \CC_{\la}, \qquad \qquad & & &
\end{align*}
where $\iota^*$ stands for pullback of functions and $\iota_*$ is
averaging along fibers, i.e. taking invariants with respect to
$\Ker(\iota)$.

\subsection{Infinitesimal induction and restriction}\label{infinitesimalinduction} The second
family of functors is an exclusive feature of higher level groups
which disappears in the case of the categories $\{ \CC_{1^n} \}$. It
stems from the fact that an $\OO$-module which is not annihilated by
$\m$ can have proper submodules of same rank.

\pagebreak

\begin{defnprop}\label{unique} Let $\mu \le \la$ be isomorphism types of finite
$\OO$-modules. The following are equivalent

\begin{enumerate}
\item (unique embedding) For every
pair of embeddings $\x \hookrightarrow \y$ and $\x' \hookrightarrow \y'$ of
modules of type $\mu$ in modules of type $\la$ and an isomorphism
$h:\x \simeq \x'$ there exists an isomorphism $\tilde{h}:\y \simeq
\y'$ which extends $h$.

\item (unique quotient) For every pair of surjections $\y
\twoheadrightarrow \x$ and $\y' \twoheadrightarrow \x'$ of modules
of type $\la$ onto modules of type $\mu$ and an isomorphism $h:\x
\simeq \x'$ there exists an isomorphism $\tilde{h}:\y \simeq \y'$
which lifts $h$.

\item ($\mathcal{G}(\mu,\y)$ is a homogeneous $\G_{\la}$-space) The group $\G_{\la}$ acts transitively on
$\mathcal{G}(\mu,\y)$, the Grassmannian of submodules of type $\mu$
in $\y$.
\end{enumerate}
In such case say that {\em $\mu$ is symmetric in $\la$}.
\end{defnprop}

\begin{proof} If $m \le n$, multiplication by $\pi^{n-m}$ is an injective $\OO$-homomorphism of $\OO_m$ into $\OO_n$.
 Let $E=\mathop{\varinjlim}\limits \OO_n$ denote the direct limit. Then $E$ is an injective $\OO$-module and for any
finite $\OO$-module $M$, its {\em dual} is defined to be
$\check{M}=\Hom_{\OO}(M,E)$. In \cite[chapter 2]{MI1} it is shown
that $M \mapsto \check{M}$ is an exact contravariant functor, that
$M \simeq \check{M}$, and that submodules of type $\mu$ and cotype
$\nu$ in $M$ are in bijective correspondence with submodules of type $\nu$
and cotype $\mu$ in $\check{M}$. It follows that any diagram of the
form
\begin{equation}\label{diagram.embedding}
\begin{matrix}
\x &\hookrightarrow& \y \\
 \uparrow \!\!\wr & &   \\
\x' &\hookrightarrow& \y'
\end{matrix}
\end{equation}
could be completed to a commutative diagram with vertical
isomorphisms if and only if any dual diagram of the form
\begin{equation}\label{diagram.epi}
\begin{matrix}
\check{M}_{\la} &\twoheadrightarrow& \check{M}_\mu \\
 & & \wr \!\! \downarrow    \\
\check{M}_{\la}' &\twoheadrightarrow & \check{M}_{\mu}'
\end{matrix}
\end{equation}
could be completed as well. This establishes the equivalence of (1)
and (2). Part (3) is a simple reformulation of (1).
\end{proof}

\begin{remark} \qquad

\begin{enumerate}
\item[(a)] The rectangular partitions of the form $\kk^n=(\kk,\ldots,\kk)$ enjoy
the property that every $\mu \le \kk^n$ is symmetric in them (in
\cite{BO1} such types are called symmetric themselves). The reader
is referred to \cite{BO1} for more details on
$\mathcal{G}(\mu,M_{\kk^n})$ and the representations of $G_{\kk^n}$
arising from it.

\item[(b)] Though the definition makes perfect sense for any length
 of the type $\mu$, we shall assume in the sequel that
 $\text{rank}(M_\mu)=\text{rank}(M_\la)$, or equivalently
 $\text{length}(\mu)=\text{length}(\la)$,
  in order to avoid overlap with geometric induction and restriction.

\end{enumerate}
\end{remark}

Let $\mu$ be symmetric in $\la$ with equal length. We define four
functors
\begin{align*}
& & \ie{\mu}{\la},~\iq{\la}{\mu}&: \CC_{\mu} \to \CC_{\la} & & &
\\ & &
                  \re{\mu}{\la}, ~\rqq{\mu}{\la}&: \CC_{\la} \to \CC_{\mu},     & & &
\end{align*}
as follows. Choose a submodule $\x$ of type $\mu$ in a module $\y$
of type $\la$. Let $P_{\x,\y} \subset \Aut_{\OO}(\y) $ be the
stabilizing group of $\x$. Then by the unique embedding property the
restriction map gives an epimorphism
\begin{equation}\label{ueinduction}
 \Aut_{\OO}(\y)  \supset
\Pp_{\x,\y} \overset{\varphi}{\to} \Aut_{\OO}(\x).
\end{equation}
Let $\xi$ be a representation of $\Aut_{\OO}(\x)$ and $\eta$ be a
representation of $\Aut_{\OO}(\y)$. Define
\begin{align*}
& & \ie{\mu}{\la}&: \CC_{\mu} \to \CC_{\la}
   & &\xi \mapsto \Ind_{\Pp_{\x,\y}}^{\Aut_{\OO}(\y)}(\varphi^{*}\xi)
& & & \\ & & \re{\mu}{\la}&: \CC_{\la} \to \CC_{\mu} & &\eta \mapsto
\varphi_{*}\Res_{\Pp_{\x,\y}}^{\Aut_{\OO}(\y)}(\eta).& & &
\end{align*}

Let $\z$ be a submodule of $\y$ with $[\y/\z]=\mu$. Let $P_{\z,\y}
\subset \Aut_{\OO}(\y) $ be the stabilizing group of $\z$. Then the
group $\Pp_{\z,\y}$ acts on the quotient $\y/\z$, and hence maps
into $\Aut_{\OO}(\y/\z)$. This map is in fact onto since $\mu$ has
the unique quotient property. Thus, similarly to \eqref{ueinduction}
we have
\begin{equation}\label{uqinduction}
 \Aut_{\OO}(\y)  \supset
\Pp_{\z,\y} \overset{\varepsilon}{\to} \Aut_{\OO}(\y/\z),
\end{equation}
and can define
\begin{align*}
& & \iq{\mu}{\la}&: \CC_{\mu} \to \CC_{\la}
   & &\xi \mapsto \Ind_{\Pp_{\z,\y}}^{\Aut_{\OO}(\y)}(\varepsilon^{*}\xi)
& & & \\ & & \rqq{\mu}{\la}&: \CC_{\la} \to \CC_{\mu} & &\eta
\mapsto \varepsilon_{*}\Res_{\Pp_{\z,\y}}^{\Aut_{\OO}(\y)}(\eta). &
& &
\end{align*}

\subsection{Cuspidality}

The point of view taken in this paper is that the geometric
induction functors defined above, though enabling the construction
of many representations of the group $\G_{\la}$, fall too short of
leaving out a manageable part of the irreducible representations to
be constructed by other means. However, together with infinitesimal
inductions, one obtains a wealth of irreducible representations,
hopefully leaving out a set which has a natural parametrization and can be
constructed in a uniform manner.

\begin{defn}\label{cuspidalitydef}
A primitive irreducible representation of $\G_{\la}$ is {\em
cuspidal} if the image of any of its twists by one dimensional
character under the geometric restriction functors $
\rs{\mu',\mu''}{\la}$ (with $\la=\mu' \cup \mu''$) and the
infinitesimal restriction functors $\re{\mu}{\la}, \rqq{\mu}{\la}$
(with $\mu$ symmetric in $\la$ with same height and length)
vanishes.
\end{defn}

\begin{remark}\qquad

\begin{enumerate}
\item[(a)] For the special case $\la=\kk^n$, a rectangular
partition, there are related notions of cuspidality for the group $\G_{\kk^n}$ which are oriented towards the construction of supercuspidal representations of
$\GL_n(F)$, hence their name. See e.g.,
Kutzko \cite[\S2]{Kutzko1}, Carayol
\cite[\S4]{carayol} ({\em tr\mygrave{e}s cuspidale}), Hill
\cite[\S4]{HG3}, Casselman \cite[p. 316]{Casselman} ({\em strongly cuspidal}) and
Prasad \cite[\S2]{PD1} ({\em very cuspidal}). Most of them, and in fact up to slight variants which will not be discussed here all of them, agree with the following definition. A representation of $\G_{\kk^n}$ is {\em strongly cuspidal} if its restriction to $\Ker\{\G_{\kk^n} \to \G_{(\kk-1)^n}\} \simeq M_n(\Fq)$ consists
of matrices with irreducible characteristic polynomial under the identification $M_n(\Fq) \to M_n(\Fq)\hat{~}$, given by $A \mapsto \psi \left(\text{trace}(A \cdot)\right)$, where $\psi : \Fq \to  \C^{\times}$ is a nontrivial character. In \cite{AOPS} it is shown that strongly cuspidal representations are cuspidal,
 with a reverse implication when $n$ is prime.

\item[(b)] Since the groups $\Pp_{\x,\y}$ and $\Pp_{\z,\y}$ above depend only on
$\mu$ and $\la$ up to conjugation we shall denote them $\Pp_{\mu
\hookrightarrow \la}$ and $\Pp_{\la \twoheadrightarrow \mu}$
(respectively) whenever chances for confusion are slim.

\end{enumerate}

\end{remark}

\pagebreak


\section{Specializing to the rank two case}\label{sec:rank2}
From this point onwards we specialize the discussion to automorphism
groups of finite $\OO$-modules of rank two. We shall now write down
concrete matrix realization of the groups under study. We denote
$\m^m_\kk$ for the image of $\m^m$ in $\OO_{\kk}$ ($m \le \kk$). Let
$\lambda=(\lx, \ly)$ be a type of length two, i.e. a partition with two
positive parts written in a descending order and let
$M_{\lambda}=\OO_{\lx}e_1 \oplus \OO_{\ly}e_2$ denote the
$\OO$-module of type $\lambda$ with basis $\{e_1,e_2\}$. The
endomorphism ring $\E_{\la}=\End_{\OO}\left(M_{\la}\right)$ consists
of elements
\[
f=\left(\begin{matrix} f_{11} & f_{12} \\ f_{21} & f_{22}
\end{matrix}\right) \in \left[\begin{matrix} \Hom_{\OO}(\OO_{\lx} e_1, \OO_{\lx} e_1) & \Hom_{\OO}(\OO_{\ly} e_2,
\OO_{\lx} e_1) \\
\Hom_{\OO}(\OO_{\lx} e_1, \OO_{\ly} e_2) & \Hom_{\OO}(\OO_{\ly} e_2,
\OO_{\ly} e_2)
\end{matrix}\right]
\]
acting by
\[
f(m_1,m_2)=\left(\begin{matrix} f_{11} & f_{12} \\ f_{21} & f_{22}
\end{matrix}\right)\left(\begin{matrix} m_1 \\ m_2
\end{matrix}\right)=\left(\begin{matrix} f_{11}(m_1)+f_{12}(m_2) \\
f_{21}(m_1)+f_{22}(m_2)
\end{matrix}\right).
\]
Writing $f_{ij}(e_j)=a_{ij}e_i$ with $a_{ij} \in \OO_{\ell_i}$,
observe that composition of maps $f \circ f'$ is regular
multiplication of the corresponding matrices $g_f \cdot g_{f'}$, and
we have
\begin{equation}\label{Ekj}
\E_{\la} \simeq \left[\begin{matrix} \OO_{\lx} & \m^{\lx-\ly}_{\lx} \\
\OO_{\ly} & \OO_{\ly}
\end{matrix}\right] \supset \left[\begin{matrix} \OO_{\lx}^{\times} & \m^{\lx-\ly}_{\lx}\\
 \OO_{\ly}  & \OO_{\ly}^{\times}
\end{matrix}\right] \simeq \G_{\la}
\end{equation}
if $\lx > \ly$, and
\begin{equation}\label{Ekk}
\E_{\kk^2} \simeq M_2(\OO_{\kk}) \supset GL_2(\OO_{\kk}) \simeq
\G_{\kk^2}
\end{equation}
if $\lx=\ly=\kk$.

\medskip

The \lq parabolic\rq~  subgroups which are of interest arise from pairs
$M_\mu \subset M_{\la}$ possessing the unique embedding property. These are classified in the next Lemma.

\begin{lemma} Let $\la=(\lx,\ly)$ be an element in $\Lambda_2$.

\begin{enumerate}

\item If $\lx=\ly$ then every $\mu \le \la$ is symmetric in $\la$.

\item If $\lx > \ly$ then $\mu=(m_1,m_2)$ is symmetric in $\la$ if and only if $\lx=m_1$.
\end{enumerate}

\end{lemma}

\begin{proof} The ring $\OO_\kk$ is self injective, therefore, any diagram of the form \eqref{diagram.embedding} with $\lambda=\kk^n$
can be completed to a commutative diagram with vertical isomorphisms. In particular (1) follows.

To prove part (2), first assume that $m_1=\lx$. Given inclusions $M_\mu \subset M_\la$, $M'_\mu \subset M'_\la$ and isomorphism $\varphi:M_\mu \to M'_\mu$, let $\{e_1,e_2\}$ be a basis of $M_\mu$ such that $\OO_{m_i} e_i \simeq \OO_{m_i}$ for $i=1,2$. As $e_2$ belongs to the $\pi^{m_2}$-torsion of $M_\la$, there
 exist $f_2 \in M_\la$ such that $e_2=\pi^{\ly-m_2}f_2$ and $\{e_1,f_2\}$ is a basis of $M_\la$. Similarly, $\{e'_1=\varphi(e_1),e'_2=\varphi(e_2)\}$ is a basis of $M'_\mu$, and there exist $f'_2 \in M'_\la$ such that $e'_2=\pi^{\ly-m_2}f'_2$ and \{$e'_1,f'_2\}$ is a basis of $M'_\la$. Then, $\tilde{\varphi}:M_\la \to M'_\la$ defined by $\tilde{\varphi}(e_1)=\varphi(e_1)=e'_1$ and $\tilde{\varphi}(f_2)=f'_2$ is the desired extension of $\varphi$.

 Conversely, if $m_1 < \lx$, we exhibit two inequivalent embeddings; let $\{f_1,f_2\}$ be a basis of $M_\la$ such that $\OO_{\kk_i} f_i \simeq \OO_{\kk_i}$ ($i=1,2$)
 and let $e_i=\pi^{\kk_i-m_i}f_i$  be a basis of a submodule of type $\mu=(m_1,m_2)$.
  We separate into two cases: if $\ly \le m_1$ the map $\varphi$ defined by $\varphi(e_1)=e_1+f_2$ and $\varphi(e_2)=e_2$ is an isomorphism of $M_\mu$ onto $\varphi(M_\mu)$ which has no extension, since $e_1 \in \pi M_\la$ but $\varphi(e_i) \in \!\!\!\!\!/ ~~\pi M_\la$; if $\ly>m_1$ the map $\varphi$ defined by $\varphi(e_1)=\pi^{\ly-m_1}f_2$ and $\varphi(e_2)=\pi^{\lx-m_2}f_1$ is an isomorphism of $M_\mu$ onto $\varphi(M_\mu)$ which has no extension, since $e_1 \in \pi^{\lx-m_1}M_\la$ but $\varphi(e_1)\in \!\!\!\!\!/ ~~\pi^{\lx-m_1} M_\la$.

\end{proof}

The case $\mu=(\lx,0)$ corresponds
to the direct sum decomposition which contributes via the geometric
induction functor. The other types with $0 < m_2 < \ly$ and $\lx=m_1$, which are symmetric in $\la$ and have the same length and height, are denoted $\PP_\la$, and contribute via the
infinitesimal induction functors. Unraveling definitions and using the matrix
representation of the groups, the infinitesimal stabilizers are
given by
\begin{align}
& &\Pp_{\mu \hookrightarrow \la} &= \left[\begin{matrix} \OO_{\lx}^{\times} & \m^{\lx-\ly}_{\lx}  \\
\m^{\ly-m_2}_{\ly} & \OO_{\ly}^{\times}
\end{matrix}\right] \twoheadrightarrow \left[\begin{matrix} \OO_{\lx}^{\times} & \m^{\lx-\ly}_{\lx-\ly+m_2} \\
\m^{\ly-m_2}_{\ly} & \OO_{m_2}^{\times}
\end{matrix}\right] \simeq \G_{\mu} & & & \label{first}\\
& &\Pp_{\la \twoheadrightarrow \mu}&= \left[\begin{matrix} \OO_{\lx}^{\times} & \m^{\lx-m_2}_{\lx} \\
\OO_{\ly} & \OO_{\ly}^{\times}
\end{matrix}\right] \twoheadrightarrow \left[\begin{matrix} \OO_{\lx}^{\times} & \m^{\lx-m_2}_{\lx} \\
\OO_{m_2} & \OO_{m_2}^{\times}
\end{matrix}\right] \simeq \G_{\mu}. & & & \label{second}
\end{align}
The subgroup $\Pp_{\mu \hookrightarrow \la}$ consists of elements in
$G_{\la}$ which stabilize $M_{\mu}=\OO_{\lx}e_1 \oplus
\m^{\ly-m_2}_{\ly} e_2$, and the epimorphism is restriction to
$M_{\mu}$. The subgroup $\Pp_{\la \twoheadrightarrow \mu}$ consists
of elements in $\G_{\la}$ which stabilize the unique embedding of a
module of type $\la/\mu=(0,\ly-m_2)$ in $M_{\la}$ (up to
automorphism) with cotype $\mu$, namely
$M_{\la/\mu}=\m_{\ly}^{m_2}e_2 \subset M_{\la}$, and the epimorphism
is the action of these elements on $M_{\la}/M_{\la/\mu}=\OO_{\lx}e_1
\oplus \OO_{m_2} e_2$.

If $m_2=0$ the groups $\Pp_{\mu \hookrightarrow \la}$ and $\Pp_{\la
\twoheadrightarrow \mu}$ become the geometric stabilizers of
$\OO_{\lx}e_1$ and $\OO_{\ly}e_2$ (respectively) in $M_{\la}$:
\begin{align}
& &\Pp_{\mu \hookrightarrow \la} &=\Pp_{(\lx),(\lx,\ly)}= \left[\begin{matrix} \OO_{\lx}^{\times} & \m^{\lx-\ly}_{\lx}  \\
 & \OO_{\ly}^{\times}
\end{matrix}\right] & & & \label{third}\\
& &\Pp_{\la \twoheadrightarrow \mu} &=\Pp_{(\ly),(\lx,\ly)}= \left[\begin{matrix} \OO_{\lx}^{\times} &  \\
\OO_{\ly} & \OO_{\ly}^{\times}
\end{matrix}\right]. & & & \label{forth}
\end{align}


\section{Representations of $ \G_{(\kk,1)}$}\label{sec:Gk1}
The group $\G_{(\kk,1)}$ is somewhat different from the groups
$\G_{(\lx,\ly)}$ with $\ly>1$. The fundamental difference from the
point of view taken here is that the geometric induction suffices to
build the induced part. We treat this group separately since its
analysis is easier than the other automorphism groups. It plays an
important role, being the first step in the inductive scheme which
is developed later on. The aim of this section is to classify the
irreducible representations of $\G_{(\kk,1)}$ and to identify the
cuspidal representations. We assume that $\kk>1$, as the case $\kk=1$ is well
known and also does not fall into the present setup. We have
\[
\begin{matrix}
& (\Oa,+)  & \hookrightarrow   &   H(\Oa) & \hookrightarrow &
\G_{(\kk,1)} & \twoheadrightarrow & \G_{(\kk-1)} \times \G_{(1)}&
\\
&    \wr\!\!\parallel  &  &\wr \!\!\parallel&  &\parallel & & \parallel& \\
Z= & \bigg[\begin{matrix} 1+\m^{\kk-1}_{\kk} &  \\
  & 1
\end{matrix}\bigg] &\hookrightarrow &  \bigg[\begin{matrix} 1+\m^{\kk-1}_{\kk} & \m^{\kk-1}_{\kk} \\
 \Oa & 1
\end{matrix}\bigg] &  \hookrightarrow & \bigg[\begin{matrix} \OO_{\kk}^{\times} & \m^{\kk-1}_{\kk} \\
 \Oa & \OO_{1}^{\times}
\end{matrix}\bigg]& \overset{p}{\twoheadrightarrow} &\OO_{\kk-1}^{\times} \times \Oas
\end{matrix}
\]
where $H=H(\Oa)$ stands for the $\Oa$-points of the Heisenberg
group. One verifies directly that indeed $\Ker(p)
\overset{\sim}{\to} H(\Oa)$ by
\[
\begin{matrix}
\bigg(\begin{matrix} 1+\pi^{\kk-1}u & \pi^{\kk-1}v \\
 w & 1
\end{matrix}\bigg) \mapsto \Bigg(\begin{matrix} 1 & v &  u \\
 & 1 & w \\
& & 1
\end{matrix}\Bigg), \quad u,v,w \in \Oa.
\end{matrix}
\]
We shall freely identify $H$ with $\Ker(p)$ and $Z(H)$ with
$Z$.

By the Stone-von Neumann Theorem, for each nontrivial character of the center
$\chi:Z \to \C^{\times}$, there is a unique irreducible
representation $\rho_{H,\chi}$ of $H$ of dimension $q=[H:Z]^{1/2}$. The representation $\rho_{H,\chi}$ is constructed by choosing any
extension of $\chi$ to a maximal abelian subgroup of $H$ and then inducing it to $H$.

The other irreducible
representations of $H$ correspond to the trivial character and hence
factor through the quotient $H/Z \simeq \Oa \oplus \Oa$.
Summarizing, the irreducible representations of $H$ are
\[
\widehat{H}= \widehat{H/Z} ~\bigsqcup~ \{ \rho_{H,\chi} ~|~ 1 \ne \chi \in \widehat{Z} \} ,
\]

By Clifford's theorem (cf.
\cite[Theorem 6.2]{Isaacs}) the restriction of any irreducible
representation of a finite group to a normal subgroup is a multiple
of a sum over a full orbit of representations. Applying the theorem
to $H \triangleleft \G_{(\kk,1)}$ we need to analyze
the action of $\G_{(\kk,1)}$ on $\hat{H}$ and then classify the
irreducible representations of $\G_{(\kk,1)}$ which lie above each
orbit.

\smallskip

{\em Case 1:} The $q$-dimensional representations of $H$. Each of
the representations $\rho_{H,\chi}$ is in fact stabilized by $\G_{(\kk,1)}$.
Indeed, the subgroup $Z$ is central in $\G_{(\kk,1)}$, hence the
action of $\G_{(\kk,1)}$ on its characters is trivial. Since each of
the representations $\rho_{H,\chi}$ is determined by $\chi$,
it must be stable under the $\G_{(\kk,1)}$-action as well. We claim
that $\rho_{H,\chi}$ can be extended to a representation of
$\G_{(\kk,1)}$. Indeed, let $T \simeq \OO_{\kk}^{\times} \times \Oas$ be the subgroup of diagonal elements in $\G_{(\kk,1)}$, and let $H_1$ be a maximal
abelian subgroup of $H$. Then $\chi$ can be extended from $Z=T \cap H_1$ to $TH_1$, and inducing the extension from $ TH_1$ to $\G_{(\kk,1)}$ gives a $q$-dimensional representation which must extend $\rho_{H,\chi}$. Thus, there are $|G_{(\kk,1)}/H|=|T/Z|=q^{\kk-2}(q-1)^2$ different extensions of $\rho_{H,\chi}$, and altogether  $q^{\kk-2}(q-1)^3$ representations of $\G_{(\kk,1)}$ counting for all $1 \ne \chi \in \widehat{Z}$.

\smallskip

{\em Case 2:} The one-dimensional representations of $H$. Since
these representations are in fact representations of the quotient
$H/Z \simeq V = \Oa \oplus \Oa$ we need to understand the action
of $G_{(\kk,1)}/Z$ on $\hat{V}$. Fix a nontrivial character
$\psi:\Oa \to \C^{\times}$, and identify $\hat{V}$ with $V$ by
$\langle (\hat{v},\hat{w}),(v,w) \rangle = \psi(\hat{v}v+\hat{w}w)$
for $ (\hat{v},\hat{w}),(v,w) \in V$. The action of $G_{(\kk,1)}/Z$
on $V$ by conjugation gives
\[
g:(v,w) \mapsto (d^{-1}av,a^{-1}dw), \qquad (v,w) \in V, \quad g =
\left(\begin{smallmatrix} a & b\pi^{\kk-1} \\ c & d
\end{smallmatrix}\right)Z \in \G_{(\kk,1)}/Z,
\]
and then, using $\langle g \cdot (\hat{v},\hat{w}),(v,w) \rangle
=\langle (\hat{v},\hat{w}),{g^{-1}} \cdot (v,w) \rangle$, we have
\[
g:(\hat{v},\hat{w}) \mapsto (a^{-1}d\hat{v},d^{-1}a\hat{w}), \qquad
(\hat{v},\hat{w}) \in \hat{V}.
\]
There are three sporadic orbits and one family of orbits parameterized by $\Oas$:
\begin{align*}
& & \text{A}~~ & & \text{B}_+ ~~& & \text{B}_- ~~& & \text{C}=\{\text{C}_{\hat{w}} ~|~ \hat{w} \in \Oas \}  & & &\\
& & [\hat{0},\hat{0}] & & [\hat{0},\Oas]& &[\Oas,\hat{0}] & &
\left\{ \Oas\cdot[1,\hat{w}] ~|~ {\hat{w} \in \Oas}\right\} & & &
\end{align*}
The representations of $\G_{(\kk,1)}$ which lie above the character
A are one dimensional and are pullbacks from characters of
$\G_{(\kk-1)}\times \G_{(1)}$ along $p$. The stabilizer in
$\G_{(\kk,1)}$ of any of the characters in $\text{B}_{\pm}$ and C is
easily seen to be the subgroup $D H$ with $D=\{d \cdot \text{Id} ~|~d \in
\OO_{\kk}^{\times}\} \subset \G_{(\kk,1)}$. As in the previous case, they can be extended to $DH$, in $[D:Z]$ many ways, each of them induce irreducibly to $\G_{(\kk,1)}$, giving altogether
\[
(|\text{B}_+|+|\text{B}_-|+|C|)[D:Z]=(q+1)(q-1)q^{\kk-2}
\]
distinct irreducible representations of $\G_{(\kk,1)}$ of dimension
$[\G_{(\kk,1)}:DH]=q-1$.

\begin{theorem} \qquad
\begin{enumerate}
\item $\RR_{(\kk,1)}(\DD)=q^{\kk-2}(q-1)^2\DD^1 + q^{\kk-2}(q^2-1)\DD^{q-1} + q^{\kk-2}(q-1)^3\DD^q$.

\item The $q^{\kk-2}(q-1)^2$ representations of dimension $(q-1)$ which lie
above the orbit C are cuspidal. All the other representations,
namely those lying above types A, $\text{B}_{\pm}$ and the
$q$-dimensional representations, are induced.
\end{enumerate}
\end{theorem}

\begin{proof} Part (1) follows from the discussion above and supplies a proof of Theorem \ref{poly} part (1) from the
introduction. As for part (2), the representations of types A and
$\text{B}_{\pm}$ contain the trivial character of either
\begin{equation}\label{Upm}
U_+=\bigg[\begin{matrix} 1 & \m^{\kk-1}_{\kk} \\
   & 1
\end{matrix}\bigg] \qquad \text{or} \qquad U_-=\bigg[\begin{matrix} 1 &  \\
 \Oa  & 1
\end{matrix}\bigg].
\end{equation}
It follows that applying one of the functors
$\rs{(\kk),(1)}{(\kk,1)}$ or $\rs{(1),(\kk)}{(\kk,1)}$ to these
representations gives a nonzero representation of $\G_{(\kk)} \times
\G_{1}$, hence they are geometrically induced. The representations
of type C do not contain the trivial character of the subgroups in
\eqref{Upm} and are hence annihilated by both functors. Since these
are the only restriction functors for the type $(\kk,1)$ we conclude
that these representations are cuspidal. As for the $q$-dimensional representations, we know that each of them is an extension of $\rho_{H,\chi}$. One of the ways to construct the latter is by inducing a character of $ZU_{\pm}$ which is trivial on $U_{\pm}$ and nontrivial on $Z$, in particular $\rho_{H,\chi}$ and any of its extensions contains the trivial representation upon restriction to $U_{\pm}$.

\end{proof}


\section{The induction scheme}\label{sec:inductionscheme}

To set up the induction scheme we start by analyzing the orbits of
$\G_{\la}$ on the reduction map kernel. For any $\la \in \Lambda_n$ which has the property that $\f{\la}$ is also in $\Lambda_n$
let
\[
K_{\la} \hookrightarrow \G_{\la} \twoheadrightarrow \G_{\f{\la}},
\]
be defined as follows. The type of $\pi M_\la$ is $\f{\la}$, and the
embedding $\pi M_\la \subset M_\la$ has the property that any
automorphism of $\pi M_\la$ can be extended to an automorphism of
$M_{\la}$. It follows that restriction maps $\G_{\la}$ onto
$\G_{\f{\la}}$. The kernel of the map is
$K_\la = \text{Id}+\Hom_\OO\left(M_\la, M_\la[\pi] \subset M_\la\right)$,
where $M_\la[\pi]$ stands for the $\pi$-torsion points in $M_{\la}$, and can be identified with $\left(M_n(\OO_1), + \right)$.
Appealing to Clifford theory, the irreducibles of $\G_{\la}$ can
hence be naturally partitioned according to the orbits they mark on
$\hat{K}_{\la}$.

Let $(2,2) \le \la \in \Lambda_2$. We shall now describe the action
of $\G_{\la}$ on $\hat{K}_{\la}$ explicitly. Using the matricial form \eqref{Ekj} the identification $K_{\la}
\overset{\sim}{\to} M_2(\OO_1)$ is given by
\[
I+\left( \begin{matrix} \pi^{\lx-1}u & \pi^{\lx-1}v \\ \pi^{\ly-1}w
& \pi^{\ly-1}z  \end{matrix} \right) \mapsto \left( \begin{matrix} u
& v
\\ w & z \end{matrix} \right).
\]
Identify the dual ${M_2(\OO_1)}\hat{~}$ with $M_2(\OO_1)$ by $x
\mapsto \psi(\text{trace}(^\text{t}x \,\cdot))$, $\psi: \OO_1 \to
\C^{\times}$ a nontrivial character. Using these identifications,
the action of $\G_{\la}$ on $\hat{K}_{\la}$ can be written
explicitly in terms of $ {M_2(\OO_1)}\hat{~}$ as
\begin{align}
& g =\left(\begin{matrix} a & b\pl \\
c & d
\end{matrix} \right) : \left(\begin{matrix} \hat{u} & \hat{v} \\
\hat{w} & \hat{z}
\end{matrix} \right) \mapsto \left(\begin{matrix} \hat{u} & \frac{d}{a} \hat{v} + \frac{c}{a}\hat{u}
\\  \frac{a}{d}\hat{w}  -\frac{b}{d}\hat{u}   \quad&
 \hat{z} - \frac{b}{a} \hat{v}  + \frac{c}{d} \hat{w}  - \frac{
 bc}{ad}\hat{u}
  \end{matrix} \right), \qquad \lx > \ly,&  & &  \label{actionGkjonK} & & \\
& g =~\left(\begin{matrix} a & b \\
c & d
\end{matrix} \right) : \left(\begin{matrix} \hat{u} & \hat{v} \\
\hat{w} & \hat{z}
\end{matrix} \right) \mapsto \left(\begin{matrix} a & b \\
c & d
\end{matrix} \right)^{-1} \left(\begin{matrix} \hat{u} & \hat{v} \\
\hat{w} & \hat{z}
\end{matrix} \right)    \left(\begin{matrix} a & b \\
c & d
\end{matrix} \right), \qquad \lx = \ly = \kk,&  & &  \label{actionGkkonK} & &
\end{align}
where $\pl=\pi^{\lx-\ly}$. The action \eqref{actionGkjonK} is a
particular instance of the action \eqref{dualaction} which we use
here without proof, postponing the detailed discussion to
\S\ref{sec:constructioncuspidal}, whereas the action
\eqref{actionGkkonK} of $G_{\kk^2}$ is simply the action of
$\G_{1^2}=GL_2(\OO_1)$ on $M_2(\OO_1)$ by conjugation via the
reduction map $\G_{\kk^2} \twoheadrightarrow \G_{1^2}$.

The first step in the induction scheme is to identify and isolate
from the discussion the non-primitive irreducible representations of
$\G_{\la}$. These representations are twists of pullbacks of
representations of $\G_{\f{\la}}$. Essentially, there are $q=|\Oa|$
such twists which are described as follows. Let \lq $\det$\rq~
denote the determinant map of $\G_{\la}$
\[
\begin{split}
\det:~ \G_{\la} &\to \OO_{\ly}^{\times} \\
\left(\begin{smallmatrix} a & b\pl \\
c & d
\end{smallmatrix} \right) &\mapsto ad-bc\pl \pmod {\pi^{\ly}}.
\end{split}
\]
Any character $\chi:\OO_{\ly}^{\times} \to \C^{\times}$ gives
rise to a one-dimensional character $\det^*\chi$ of $\G_{\la}$.
However, from the point of view of twisting, the interest lies in
characters of $\OO_{\ly}^{\times}$ modulo the characters of
$\OO_{\ly-1}^{\times}$. These are represented by an arbitrary
extension of each of the characters
$\psi(\hat{z}\,\cdot):1+\m^{\ly-1}_{\ly} \simeq \Oa \to \C^{\times}$
($\hat{z} \in \Oa$), from $1+\m^{\ly-1}_{\ly}$ to
$\OO_{\ly}^{\times}$. Denote these representatives by
$\chi_{\hat{z}}$ ($\hat{z} \in \Oa$).

\subsection{Non-rectangular case}
Let $\la \in \Lambda_2$ be a partition with unequal parts greater than $1$.
 Using \eqref{actionGkjonK}, representatives for the orbits of $\G_{\la}$ on $\hat{K}_{\la}$, their number, and the total number of elements
 they contain, are given in the following list\footnote{The
\lq$\hb$\rq ~in the matrices' entries stands for $x \mapsto \psi(1
\cdot x)$ and is a convenient noncanonical choice of a nontrivial
character of $\OO_1$. \lq$\ha$\rq ~stand for $x \mapsto \psi(0\cdot
x) \equiv 1$, the trivial character of $\Oa$.}

\smallskip

 \centerline{\bf \underline{Table 1: Orbits in the
non-rectangular case}}
\begin{align*}
&  & &\text{\em{Representatives for $\rho_{|K}$ }} & & \text{\em{Number of orbits}} & & \text{\em{Total number of
elements}}&\\
& (i)  & &\left(\begin{matrix} \ha & \ha \\
\ha & \hat{z}
\end{matrix} \right), ~\hat{z} \in \OO_1 & & \qquad q & & \qquad q \\
&(ii) &   &\left(\begin{matrix} \hat{u} & \ha \\
\ha & \hat{z}
\end{matrix} \right), ~\hat{z} \in \OO_1, \hat{u} \in \Oas & & \quad q(q-1) & &\quad q^3(q-1) \\
&(iii) &   &\left(\begin{matrix} \ha & \ha \\
\hb & \ha \end{matrix} \right) & & \qquad 1 & & \quad q(q-1) \\
& (iv) &   &\left(\begin{matrix} \ha & \hb \\
\ha & \ha \end{matrix} \right) & & \qquad 1 & & \quad q(q-1) \\
& (v)  &   &\left(\begin{matrix} \ha & \hb \\
 \hat{w} & \ha \end{matrix} \right), ~\hat{w} \in \Oas & & ~\quad
q-1 & & \quad q(q-1)^2
\end{align*}

\smallskip

Define the following subgroups of $\G_{\la}$
\[
U_{+}=\left[\begin{matrix} 1 &  \m^{\lx-\ly}_{\lx}      \\
 &  1
\end{matrix} \right]  \qquad \text{and}\qquad
U_{-}=\left[\begin{matrix} 1 & \\
 \OO_{\ly} &  1
\end{matrix} \right],
\]
which play the role of \lq upper/lower unipotent\rq~  subgroups.
Note that they are both isomorphic to $(\OO_{\ly},+)$ but not
conjugate in $\G_{\la}$ since $\la$ is non-rectangular. Let
$V_{\pm}=U_{\pm} \cap K_{\la}$ and let $V_1 < K_{\la}$ and $V_2 <
K_{\la} $ denote the embeddings of $1+\m^{\lx-1}_{\lx}$ and
$1+\m^{\ly-1}_{\ly}$ in the $(1,1)$ and $(2,2)$ entries
respectively. Thus, $K_{\la}=V_1  V_2  V_+  V_-$.

\begin{prop}\label{tableGkj} Let $\rho$ be an irreducible
representation of $G_\la$, then
\begin{enumerate}
\item $\rho$ is a twist of a pullback from $\hat{\G}_{\f{\la}}$ $~\Longleftrightarrow~$ $\langle
\rho_{|(V_1V_-V_+)},\ha_{(V_1V_-V_+) } \rangle \ne 0 $
 ($\rho_{|K_{\la}}$ is of type $(i)$).

\item $\rho$ is a cuspidal $~\Longleftrightarrow~$
$\langle \rho_{|V_{-}},\ha_{V_{-}} \rangle = \langle
\rho_{|V_{+}},\ha_{V_{+} } \rangle = 0 \ne \langle
\rho_{|V_{1}},\ha_{V_{1}} \rangle$  ($\rho_{|K_{\la}}$ is of type
$(v)$).

\item $\rho$ is induced (up to twist) $~\Longleftrightarrow~$ $\rho_{|K_{\la}}$ is of
type $(ii)$, $(iii)$ or $(iv)$.

\end{enumerate}

\end{prop}

\pagebreak

\begin{proof} \quad

\begin{enumerate}
\item
We observe that $\det_{|K_{\la}}$ maps $V_2$
 isomorphically on $1+\m^{\ly-1}_{\ly}$, while $V_1V_+V_- = \Ker(\det)$.
We conclude that the effect of twisting a representation $\rho$ of
$\G_{\la}$ with a one dimensional character $\det^*\chi_{\hat{z}}$
is expressed on $\rho_{|K_{\la}}$ by an additive shift of the
$(2,2)$-entry. It follows that the orbits of type $(i)$ in Table 1
occur precisely for twists of pullbacks of representations of
$\G_{\f{\la}}$.

\item  Let $\rho \in \hat{C}_{\la}$. In order that the image of $\rho$ will be zero
under each of the restriction functors it must satisfy $\langle
\rho_{|V_+},\ha_{V_+} \rangle = 0 = \langle \rho_{|V_-},\ha_{V_-}
\rangle$ (use the explicit form in \S\ref{sec:rank2}). This implies,
by using the action \eqref{actionGkjonK}, that $\hat{u}=0$ and
$\hat{w},\hat{v} \in \Oas$. Then, we may choose $\hat{v}=\hb$ and
$\hat{z}=0$. Conversely, we shall see in
\S\ref{sec:constructioncuspidal} how each of the representatives
\begin{equation}\label{cusp1}
\left(\begin{matrix} \ha & \hb \\
\hat{w} & \ha
\end{matrix} \right), \qquad \hat{w} \in \OO_1^{\times}
\end{equation}
occur as the restriction of a cuspidal representation.

\item The character $\big(\begin{smallmatrix} \ha & \ha \\
\hb & \ha \end{smallmatrix} \big)$ is trivial on $V_+V_2$ and is
thus not annihilated by the functor $\re{(\lx,\ly-1)}{(\lx,\ly)}$.
This means that any representation of $\G_{(\lx,\ly)}$ lying above
orbit of type $(iii)$ is contained in an infinitesimally induced
representation of $\G_{(\lx,\ly-1)}$ using
$\ie{(\lx,\ly-1)}{(\lx,\ly)}$. Similarly, representations lying
above type $(iv)$ are contained in an infinitesimally induced
representation of $\G_{(\lx,\ly-1)}$ using
$\iq{(\lx,\ly-1)}{(\lx,\ly)}$. Same argument holds for
representations lying above orbits of type $(ii)$ after a possible
twist.

\end{enumerate}
\end{proof}

\begin{remark}
The third part of Proposition \ref{tableGkj} is given here in order
to motivate and explain the strategy. We shall later on obtain more
precise results as in Theorem \ref{maintheorem}. In particular, all
the representations of type $(ii)$ are geometrically induced without
the requirement of a twist.
\end{remark}

\subsection{Rectangular case}\label{subsec:rectangular}  Let $\la=\kk^2$ with $\kk > 1$. The orbits of the
action \eqref{actionGkkonK} are well known, given by

\centerline{\bf \underline{Table 2: Orbits in the rectangular case}}
\vspace{-0.5cm}
\begin{align*}
&  & &\text{\em{Representatives for $\rho_{|K}$ }} & & \text{\em{Number
of orbits}} & & \text{\em{Total number of
elements}}& \\
& (i) & &\left(\begin{matrix} \hat{u} & \ha \\
\ha & \hat{u}
\end{matrix} \right), ~\hat{u} \in \OO_1 & & \qquad q & & \qquad q \\
& (ii)  & &\left(\begin{matrix} \hat{u} & \ha \\
\ha & \hat{z}
\end{matrix} \right), ~{S_2}\backslash \{\hat{u},\hat{z} \in \OO_1 | \hat{u} \ne \hat{z}\} & & \quad \gfrac{1}{2}q(q-1) & &\quad \gfrac{1}{2}
q^2(q^2-1) \\
& (iii)   & &\left(\begin{matrix} \hat{u} & \hb \\
\ha & \hat{u} \end{matrix} \right),~\hat{u} \in \OO_1 & & \qquad q & & \quad  q(q^2-1) \\
& (iv) & &\left(\begin{matrix} \hat{u} & \hat{v} \\
\hb & \hat{u} \end{matrix} \right), \hat{u} \in \Oa, \hat{v} \in
\Oas
 \smallsetminus (\Oas)^2 & & ~\quad \gfrac{1}{2}q(q-1) & & \quad
\gfrac{1}{2}q^2(q-1)^2
\end{align*}

In complete analogy with Proposition \ref{tableGkj}, representations
lying over orbits of type $(i)$ are twists of pullbacks from lower
level. Representations lying over type $(iv)$ are cuspidal. The main
differences comparing to the non-rectangular case are the following.
\begin{enumerate}
\item A twist with a primitive character manifests itself on
$\hat{K}_{\kk^2}$ as a shift with a scalar matrix. As a result,
primitive representations of $\G_{\kk^2}$ are infinitesimally
induced from cuspidal representations of $\G_{\mu}$ ($\mu \in
\PP_{\kk^2}$) only up to twist.

\item There is an extra symmetry in the infinitesimal induction
which results in equivalences between $\ie{\mu}{\kk^2}$ and
$\iq{\mu}{\kk^2}$, which is simply denoted $\ib_{\mu}^{\kk^2}$ in
this case. In terms of orbits, this can be seen in the difference
between the orbit types: types $(iii)-(iv)$ in Table 1 degenerate to
one type, $(iii)$, in Table 2.
\end{enumerate}


\section{Construction of the cuspidal representations}\label{sec:constructioncuspidal}

\subsection{Non-rectangular case}
In this section the cuspidal representations of the groups
$\G_{\la}$ are constructed for $\la=(\lx,\ly)$ with $\lx > \ly > 1$.
Let
\[
\K_{\la}^{i,\sigma}=I+\left[\begin{matrix} \m^{\lx-i}_{\lx} &  \m^{\lx-i}_{\lx}\\
\m^{\ly-i+\sigma}_{\ly}  & \m^{\ly-i+\sigma}_{\ly}
\end{matrix} \right], \qquad 1 \le i \le \ly,~ \sigma \in \{0,1\}.
\]
\begin{lemma}\label{normalabelian} $\K_{\la}^{i,\sig}$ is a normal subgroup of $\G_{\la}$. If
$i  \le  (\ly+\sigma)/2 $ then $\K_{\la}^{i,\sig}$ is abelian, and
\begin{equation*}\label{Alambda}
\K_{\la}^{i,\sigma} \simeq L_{i,\sigma} = \left(  \left[\begin{matrix} \OO_{i} &  \OO_{i} \\
\OO_{i-\sig} & \OO_{i-\sig}
\end{matrix} \right] , + \right), \qquad I+\left(\begin{matrix} \pi^{\lx-i} u & \pi^{\lx-i} v \\
\pi^{\ly-i+\sig}w & \pi^{\ly-i+\sig} z
\end{matrix} \right) \mapsto \left(\begin{matrix} u & v \\
w &  z
\end{matrix} \right).
\end{equation*}
\end{lemma}

\begin{proof}
Straightforward.
\end{proof}
{\noindent} We shall proceed in the following steps.

\begin{enumerate}
\item  Analysis of the action of $\G_{\la}$ on
$\hat{K}_{\la}^{i,\sig}$.

\item  Identification of the possible orbits of characters of
a \lq large\rq~ normal abelian subgroup ${K}_{\la}^{\ls,\epsilon}$ which can give rise to cuspidal representations.

\item  Construction of the cuspidal representations by extending these characters from
$K_{\la}^{l,\epsilon}$ to their normalizer and then inducing
them to $G_{\la}$.

\end{enumerate}

\medskip

{\bf Step 1.} Assume from now on that $i \le (\ly+\sig)/2$, so that
lemma \ref{normalabelian} is satisfied. The action of $\G_{\la}$ on
$\K_{\la}^{i,\sig}$ by conjugation translates by the above
identification to an action on $L_{i,\sig}$, explicitly given
 by
\begin{equation*}\label{action}
g : \left(\begin{matrix} u & v \\
w & z
\end{matrix} \right) \mapsto e^{-1}\left(\begin{matrix} u- \frac{c}{d} v + \frac{b\pll}{a} w - \frac{bc\pll}{ad}z
& \frac{a}{d}v + \frac{b\pll}{d}z -\frac{b\pl}{d}u
-\frac{b^2\pl\pll}{ad} w
\\ \\
\frac{d}{a} w - \frac{c}{a}z +\frac{c\pl}{a\pll}u -
\frac{c^2\pl}{ad\pll}v  \quad&
 z - \frac{b\pl}{a} w  + \frac{c\pl}{d\pll} v  - \frac{
 bc\pl^2}{ad\pll}u
  \end{matrix} \right),
\end{equation*}
where $g=\left(\begin{smallmatrix} a & b\pl \\
c  & d
\end{smallmatrix} \right)$, $e=1-a^{-1}d^{-1}bc\pl$ and
$\pl=\pi^{\lx-\ly}$. Note that the inverse of $g$ is given by
\[
g^{-1}=e^{-1}\left(\begin{matrix} a^{-1} & -a^{-1}d^{-1}b\pl \\
-a^{-1}d^{-1}c  & d^{-1}
\end{matrix} \right).
\]

In the sequel we shall need to analyze the action of $\G_{\la}$ also
on the characters of $L_{i,\sig}$. Let $\psi:\OO_{i} \to
\C^{\times}$ be a fixed primitive additive character. We identify $
L_{i,\sig}$ with its dual by
\[
\langle {\left(\begin{matrix} \hat{u} & \hat{v} \\
\hat{w} & \hat{z}
\end{matrix} \right)} , \left(\begin{matrix} u & v \\
w & z
\end{matrix} \right) \rangle = \psi\left(\hat{u}u+\hat{v}v+\pll(\hat{w}w + \hat{z}z)  \right).
\]
The action of $g \in \G_{\la}$ on the dual is defined by $ \langle g
\cdot {\theta},t \rangle= \langle {\theta},g^{-1} \cdot t \rangle $,
and is given explicitly by
\begin{equation}\label{dualaction}
g : \left(\begin{matrix} \hat{u} & \hat{v} \\
\hat{w} & \hat{z}
\end{matrix} \right) \mapsto e^{-1}\left(\begin{matrix} \hat{u}+ \frac{b\pl}{a} \hat{v} - \frac{c\pl}{d} \hat{w} -
\frac{bc\pl^2}{ad}\hat{z} & \frac{d}{a} \hat{v} -
\frac{c\pl}{a}\hat{z} +\frac{c}{a}\hat{u} - \frac{c^2\pl}{ad}\hat{w}
\\ \\
 \frac{a}{d}\hat{w} + \frac{b\pl}{d}\hat{z} -\frac{b}{d}\hat{u}    -\frac{b^2\pl}{ad} \hat{v} \quad&
 \hat{z} - \frac{b}{a} \hat{v}  + \frac{c}{d} \hat{w}  - \frac{
 bc}{ad}\hat{u}
  \end{matrix} \right).
\end{equation}

Equipped with the action of $\G_{\la}$ on the various groups
$\hat{\K}_{\la}^{i,\sig}$ we proceed to the construction of the
cuspidal representations in several steps. We repeatedly make use of
Clifford theorem with respect to the normal subgroups
$\K_{\la}^{i,\sig}$.

\medskip

{\bf Step 2.} Let $\epsilon=\ly \pmod 2$ be the parity of $\ly$ and
let $\ls=(\ly+\epsilon)/2$. We wish to extend characters of the type
\eqref{cusp1} from $K_{\la}^{1,0}$ to $K_{\la}^{\ls,\epsilon}$
($K_{\la}^{1,0}$ is denoted by $K$ in \S\ref{sec:inductionscheme}).
Let $\eta_{\hat{u},\hat{w}}$ be the following characters of
$L_{\ls,\epsilon}$
\begin{equation}\label{cuspl}
\left(\begin{matrix} \hat{u} & \hat{1} \\
 \hat{w} &  \hat{0}
\end{matrix} \right), \qquad \hat{w} \in \OO_{\ls-\epsilon}^{\times}, ~\hat{u} \in \m_{\ls}.
\end{equation}

\begin{prop} The characters $\eta_{\hat{u},\hat{w}}$ ($\hat{w} \in \OO_{\ls-\epsilon}^{\times}, ~\hat{u} \in \m_{\ls}$)
form
 a complete set of representatives for the orbits of $\G_{\la}$ on characters of $\hat{L}_{\ls,\epsilon}$ which
lie above the characters \eqref{cusp1} of $\hat{L}_{1,0}$.
\end{prop}

\begin{proof} Define the following functions on $\hat{L}_{\ls,\epsilon}$
\[
\begin{split}
\text{Tr}_{\pl}\left( \begin{matrix} \hat{u} & \hat{v} \\
\hat{w} & \hat{z}
\end{matrix}  \right)&=\hat{u}+\pl \hat{z} \pmod
{\pi^{\ls}}  \\
\text{Det}\left( \begin{matrix} \hat{u} & \hat{v} \\
\hat{w} & \hat{z}
\end{matrix}  \right)&=\hat{u}\hat{z}-\hat{w}\hat{v} \pmod
{\pi^{\ls-\epsilon}}.
\end{split}
\]
Using the action \eqref{dualaction} one easily checks that
$\text{Tr}_\pl$ and $\text{Det}$ are invariants of the
$\G_{\la}$-orbits. It follows that the characters $\eta_{\hat{u},\hat{w}}$ lie in
distinct $\G_{\la}$-orbits since
$\text{Tr}_{\pl}(\eta_{\hat{u},\hat{w}})=\hat{u}$ and
$\text{Det}_{\pl}(\eta_{\hat{u},\hat{w}})=-\hat{w}$. Conversely, one
easily verifies that any element which lie above the characters
\eqref{cusp1} can be brought to an element of the form
$\eta_{\hat{u},\hat{w}}$ the assertion follows.

\end{proof}

{\bf Step 3.} Let $N_{\la}=N_{\la}(\hat{u},\hat{w})=
\Stab_{\G_{\la}}(\eta_{\hat{u},\hat{w}})$ be the stabilizer of
$\eta_{\hat{u},\hat{w}}$ in $\G_{\la}$. The following theorem
completes the construction of the cuspidal representations of
$\G_{\la}$.

\begin{theorem}\label{cuspidaltheorem} \qquad
\begin{enumerate}

\item The normalizer of $\eta_{\hat{u},\hat{w}}$ in $\G_{\la}$ is
given by
\[
N_{\la}(\hat{u},\hat{w})=\left\{ \left(\begin{matrix} a & b\pl \\
 c &  d
\end{matrix} \right) ~\big|~ b \equiv c\hat{w} \pmod {\pi^{\ls-\epsilon}}, ~d \equiv a-c\hat{u} \pmod {\pi^{\ls}} \right\}.
\]

\item The character $\eta_{\hat{u},\hat{w}}$ can be extended to
$N_{\la}(\hat{u},\hat{w})$.

\item For any extension $\tilde{\eta}$ of $\eta_{\hat{u},\hat{w}}$ the induced representation
$\rho_{\tilde{\eta}}=\Ind_{N_{\la}}^{\G_{\la}}\left(\tilde{\eta}_{\hat{u},\hat{w}}\right)$
is irreducible and cuspidal and distinct $\tilde{\eta}$'s give rise
to inequivalent representations. Any cuspidal representation of
$\G_{\la}$ is of this form. There exist $q^{\lx+\ly-3}(q-1)^2$ such
representations and their dimension is $q^{\ly-1}(q-1)$.
\end{enumerate}

\end{theorem}

\begin{proof}\qquad

\begin{enumerate}
\item Using the action \eqref{dualaction} we see that if $g \cdot
\eta_{\hat{u},\hat{w}} = \eta_{\hat{u},\hat{w}}$, we get
\[
\begin{split}
& \qquad \text{(1,1)} \qquad \hat{u}=\frac{\hat{u}}{e}+
\frac{b\pl}{ae} \hat{1} - \frac{c\pl}{de}\hat{w}  \qquad \text{(in
$\OO_{\ls}$)} \\
& \qquad \text{(1,2)} \qquad
\hat{1}=\frac{d}{ae}\hat{1}+\frac{c}{ae}\hat{u} -\frac{c^2\pl}{ead}
\hat{w} \qquad \text{(in $\OO_{\ls}$)} \\
& \qquad \text{(2,2)}
\qquad \hat{0}= - \frac{b}{a} \hat{1}  + \frac{c}{d}\hat{w} - \frac{
 bc}{ad}\hat{u} \qquad \text{(in
$\OO_{\ls-\epsilon}$)}.
\end{split}
\]

Adding up $[(2,2)-d^{-1}be \cdot \text{(1,2)}]$ gives
\begin{equation}\label{cequiv}
 b \equiv c\hat{w} \pmod {\pi^{\ls-\epsilon}}.
\end{equation}
Equation (1,1) is equivalent to
\[
\frac{c\pl}{d}\hat{w}= \frac{b\pl}{a}\hat{1}
+\frac{bc\pl}{ad}\hat{u}, \qquad \text{(in $\OO_l$)}
\]
and the latter gives after substitution in (2,1)
\begin{equation}\label{dequiv}
d \equiv a-c\hat{u} \pmod {\pi^{\ls}}.
\end{equation}
Conversely, if relations \eqref{cequiv} and \eqref{dequiv} hold then
$g \cdot \eta_{\hat{u},\hat{w}} = \eta_{\hat{u},\hat{w}}$.

\item Let
\[
A_\la=A_{\la}(\hat{u}',\hat{w}')=\left\{\left(\begin{matrix} a & c\hat{w}'\pl \\
c  & a-c\hat{u}'
\end{matrix} \right)~\mid~ a \in \OO_{\lx}^{\times}, c \in \OO_{\ly}  \right\}.
\]
where $\hat{u}', \hat{w}' \in \OO_{\ly}$ are lifts of
$\hat{u},\hat{w}$ respectively. We have the following diagram
\[
\begin{matrix} & & G_{\la} & & \\
& & \mid & & \\
& & K_{\la}^{\ls,\epsilon}A_{\la} & & \\
& \diagup & & \diagdown & \\
K_{\la}^{\ls,\epsilon} & & & & A_{\la} \\
& \diagdown & & \diagup & \\
& & K_{\la}^{\ls,\epsilon} \cap A_{\la} & &\\
\end{matrix}
\]
Note that $N_{\la}=K_{\la}^{\ls,\epsilon}A_{\la}$, and since
$A_{\la}$ is abelian, we can extend $\eta_{\hat{u},\hat{w}}$ to $N_{\la}$ in
$[A_{\la}:K_{\la}^{\ls,\epsilon} \cap
A_{\la}]=[K_{\la}^{\ls,\epsilon}A_{\la} : K_{\la}^{\ls,\epsilon}]$
different ways.

\item Going over all the possible
choices of $\hat{u}$ and $\hat{w}$ and all the possible extensions
to the stabilizers $N_{\la}(\hat{u},\hat{w})$, we get
\[
|\m_{\ls}||\OO_{\ls-\epsilon}^{\times}||A_{\la}/K_{\la}^{\ls,\epsilon}
\cap A_{\la}|=q^{\lx+\ly-3}(q-1)^2
\]
different characters. By \cite[Theorem 6.11]{Isaacs} all these
characters induce irreducibly to $\G_{\la}$. Since all the possible
candidates of orbits in $\hat{K}_{\la}^{\ls,\epsilon}$ which might lie below a cuspidal representation of $\G_\la$ are covered, these
induced representations exhaust the cuspidal representations of
$\G_{\la}$. Their dimension is $|\G_{\la}/N_{\la}|=q^{\ly-1}(q-1)$.

\end{enumerate}

\end{proof}

\subsection{Cuspidal representations of $\G_{\kk^2}$} The cuspidal
representations of $\G_{\kk^2}$ were constructed by Nobs in
\cite[\S4]{Nobs1}  using the Weil representation under the name {\em
la s\myacute{e}rie non ramifi\myacute{e}e} (see also the more general treatment in \cite{Gerardin}). There are $\gfrac{1}{2}(q^2-1)(q-1)q^{2\kk-3}$ such representations and their dimension is
$q^{\kk-1}(q-1)$. The number of these representations matches
the number of Galois orbits of strongly primitive characters of
$\tilde{\OO}_{\kk}^{\times}$ (the latter is the units in the
non-ramified extension of $\OO_\kk$ of degree $2$). This is of no
coincidence. In fact, there exist a canonical bijective
correspondence between such orbits of characters and the cuspidal
representations of $\G_{\kk^2}$. This is proved in general in \cite{AOPS} for
types $\la=\kk^n$ where $n$ is prime. The reader is referred to {\em loc.
cit.} for another construction of the cuspidal representation of
$\G_{\kk^2}$ which is more in the spirit of this paper and avoids
the use of the Weil representation. In short, the correspondence is
established as follows. A character of
$\theta:\tilde{\OO}^{\times}_{\kk} \to \C^{\times}$ is {\em strongly
primitive} if its restriction to $\Ker\{\OO_{\kk}^{\times} \to
\OO_{\kk-1}^{\times}\} \simeq (\tilde{\OO}_1,+)$ gives a primitive
additive character of $\tilde{\OO}_1$. Using the natural embedding
$\tilde{\OO}^{\times}_{\kk} \hookrightarrow \G_{\kk^2}$, it is shown
that the centralizer of a strongly primitive character $\theta$ is
$N_{\lfloor \kk/2 \rfloor} \cdot \tilde{\OO}^{\times}_{\kk}$, where
$N_i=\Ker\{\G_{\kk^2} \to \G_{i^2}\}$, and that there exists a unique
representation $\rho_{\theta}$ of the centralizer which restricts to
$\theta$ or to one of its Galois conjugates. Inducing the
$\rho_{\theta}$'s to $\G_{\kk^2}$ gives the cuspidal
representations.


\section{Infinitesimal and geometric induction}\label{sec:inductions}

In this section we take a closer look at the infinitesimal and
geometric inductions. The first subsection describes some basic
properties of geometric and infinitesimal inductions in the general
setting. We then specialize to the rank two case and in particular
obtain finer results regarding the cases in which such inductions
are irreducible. Throughout we use the notation
$K_{\la}=\Ker\{\G_{\la} \to \G_{\f{\la}}\}$.

\subsection{Basic properties}  All the induction functors defined in
\S\ref{sec:operations} are compositions of inflation and usual
induction: $\Ind \circ \Inf$. Similarly, all the restriction
functors are compositions of usual restriction and a functor of
fixed points: $\Inv \circ \Res$. Since $(\Ind,\Res)$ form an adjoint
pair, and so do $(\text{Inf},\text{Inv})$, it follows that
\[
(\is{\mu,\mu'}{\la},\rs{\mu,\mu'}{\la}), ~
(\ie{\mu}{\la},\re{\mu}{\la}), ~\text{and}~
(\iq{\mu}{\la},\rqq{\mu}{\la}),
\]
are adjoint pairs as well. In order to prove some associativity
properties of these functors we shall need the following lemma.

\begin{lemma}\label{inf-ind} Let $G$ be a finite group and $G \supset P_1 \supset P_2
\supset  U_2 \supset U_1$ a chain of subgroups with $U_i
\triangleleft P_i$. Then
\[
\Ind_{P_2}^{G} \circ \Inf_{P_2/U_2}^{P_2}=\Ind_{P_1}^{G}
 \circ
\Inf_{P_1/U_1}^{P_1} \circ \Ind_{P_2/U_1}^{P_1/U_1} \circ
\Inf_{P_2/U_2}^{P_2/U_1}.
\]

\end{lemma}

\begin{proof} By \cite[Lemma 3.4]{Yoshida} we have $\Inf_{P_1/U_1}^{P_1} \circ \Ind_{P_2/U_1}^{P_1/U_1} =
\Ind_{P_2}^{P_1} \circ \Inf_{P_2/U_1}^{P_2}$, which together with
the associativity of induction and inflation proves the assertion.
\end{proof}

\begin{claim}\label{associativity-infinitesimal} Let $\nu \le \mu \le \la$ be types in $\Lambda$ with $\nu$
symmetric in $\mu$ and $\mu$ symmetric in $\la$. Assume that any module
of type $\nu$ is contained in a unique module of type $\mu$. Then
\begin{enumerate}

\item $\ie{\mu}{\la}  \circ \ie{\nu}{\mu}= \ie{\nu}{\la}$. Dually, $\re{\nu}{\mu} \circ \re{\mu}{\la} =
\re{\nu}{\la}$.

\item $\iq{\mu}{\la} \circ \iq{\nu}{\mu} = \iq{\nu}{\la}$. Dually, $\iq{\nu}{\mu} \circ \iq{\mu}{\la} =
\iq{\nu}{\la}$.

\end{enumerate}
\end{claim}

\begin{proof} Let $M_{\nu} \subset M_{\mu} \subset M_{\la}$ be of types $\nu$,
$\mu$ and $\la$. We should show that
\[
\Ind_{\Pp_{M_{\nu},M_{\la}}}^{\G_{\la}} \circ
\Inf_{\G_{\nu}}^{\Pp_{M_{\nu},M_{\la}}}=\Ind_{\Pp_{M_{\mu},M_{\la}}}^{\G_{\la}}
 \circ
\Inf_{\G_{\mu}}^{\Pp_{M_{\mu},M_{\la}}} \circ
\Ind_{\Pp_{M_{\nu},M_{\mu}}}^{\G_{\mu}}  \circ
\Inf_{\G_{\nu}}^{\Pp_{M_{\nu},M_{\mu}}}.
\]
By assumption we have the inclusion $\Pp_{M_{\nu},M_{\la}} \subset
\Pp_{M_{\mu},M_{\la}}$, and
\[
\begin{matrix}
U_{M_{\nu},M_{\la}} & \hookrightarrow & P_{M_{\nu},M_{\la}} &
\twoheadrightarrow & \G_{\nu} \\
\cup & & \cap & &\\
U_{M_{\mu},M_{\la}} & \hookrightarrow & P_{M_{\mu},M_{\la}} &
\twoheadrightarrow & \G_{\mu}.
\end{matrix}
\]
The first part of assertion (1) follows now from Lemma \ref{inf-ind}
using the identifications
\[
\begin{split}
\Pp_{M_{\nu},M_{\mu}}&=\Pp_{M_{\nu},M_{\la}}/U_{M_{\mu},M_{\la}} \\
\G_{\mu}&=\Pp_{M_{\mu},M_{\la}}/U_{M_{\mu},M_{\la}}\\
\G_{\nu}&=\Pp_{M_{\nu},M_{\la}}/U_{M_{\nu},M_{\la}}=\Pp_{M_{\nu},M_{\mu}}/U_{M_{\nu},M_{\mu}},
\end{split}
\]
and the second part by uniqueness of the adjoint. Part (2) follows
along the same lines by replacing the modules $M_{\nu} \subset
M_{\mu}$ by modules $M_{\la/\mu} \subset M_{\la/\nu}$ of cotypes
$\mu$ and $\nu$ respectively.
\end{proof}

\subsection{Infinitesimal induction in rank two}

\begin{prop}\label{dimensionnilpotent} Let $\la=(\lx,\ly) \in
\Lambda_2$ and let $\rho$ be an irreducible representation of
$\G_{\la}$.
\begin{enumerate}

\item If $\lx>\ly >1$ and either $\rho_{|K_{\la}} \geqslant \big(\begin{smallmatrix} \ha & \ha
\\ \hb & \ha
\end{smallmatrix}\big)$ or $\rho_{|K_{\la}} \geqslant \big(\begin{smallmatrix} \ha &
\hb
\\ \ha & \ha
\end{smallmatrix}\big)$, then $\dim(\rho)=q^{\ly-1}(q-1)$.

\item If $\lx=\ly=\kk >1$ and $\rho_{|K_{\kk^2}} \geqslant \big(\begin{smallmatrix} \ha & \hb
\\ \ha & \ha
\end{smallmatrix}\big)$, then $\dim(\rho)=q^{\kk-2}(q^2-1)$.
\end{enumerate}

\end{prop}

\begin{proof}\quad

\begin{enumerate}

\item We argue by induction on $\ly$. Let $\rho$ be an irreducible
representation of $\G_{\la}$ which contains the character
$\big(\begin{smallmatrix} \ha & \ha
\\ \hb & \ha
\end{smallmatrix}\big) \in \hat{K}_{\la}$. Using the action \eqref{actionGkjonK} we verify that this orbit
contains the $q(q-1)$ characters $\big\{\big(\begin{smallmatrix} \ha
& \ha
\\ a\hb & c\hb
\end{smallmatrix}\big)~|~ a \in \Oas,~c \in \Oa \big\}$,
out of which $(q-1)$ survive the application of the functor
$\re{(\lx,\ly-1)}{(\lx,\ly)}$, namely those having $c=0$. It follows
that if $\xi= \re{(\lx,\ly-1)}{(\lx,\ly)}(\rho)$ then $\dim
(\xi)={\dim(\rho)}/{q}$. Using the explicit form \eqref{first} we
observe that $\langle \xi_{|V_-}, \ha_{V_-} \rangle = 0$, which
means that $\xi_{|K_{(\lx,\ly-1)}}$ contains the orbit of
$\big(\begin{smallmatrix} \ha & \ha
\\ \hb & \ha
\end{smallmatrix}\big)$ or a cuspidal orbit (that is, orbits $(iii)$ or $(iv)$ in Table 1). In either
case, using the induction hypothesis or the dimension of cuspidal
representations, we get that $\xi$ is irreducible and that
$\dim(\xi)=q^{\ly-2}(q-1)$, from which the assertion follows. The
basis of the induction is $\ly=2$ in which case the role of
$K_{(\lx,1)}$ is taken over by $H$ of \S\ref{sec:Gk1} and $\xi$ lies
above the orbits $\text{B}_{+}$ or C. If $\rho$ lies above the orbit
of $\big(\begin{smallmatrix} \ha & \hb
\\ \ha & \ha
\end{smallmatrix}\big)$ we argue similarly using the
functor $\rqq{(\lx,\ly-1)}{(\lx,\ly)}$.

\item Assume $\rho \in \hat{\G}_{\kk^2}$ such that $\rho_{|K_{\kk^2}} \geqslant \big(\begin{smallmatrix} \ha & \hb
\\ \ha & \ha
\end{smallmatrix}\big)$. This orbit contains the $(q^2-1)$
characters $\big\{\big(\begin{smallmatrix} ab\hb & -a^2\hb
\\ b^2\hb & -ab\hb
\end{smallmatrix}\big)~|~ a,b \in \Oa, (a,b) \ne (0,0) \big\}$
(using the action \eqref{actionGkkonK}), out of which $(q-1)$
survive the application of the functor
$\re{(\kk,\kk-1)}{(\kk,\kk)}$, namely those having $a=0$. It follows
that if $\xi= \re{(\kk,\kk-1)}{(\kk,\kk)}(\rho)$ then $\dim
(\xi)={\dim(\rho)}/(q+1)$. Observing that $\xi_{|K_{(\kk,\kk-1)}}
\geqslant \big(\begin{smallmatrix} \ha & \ha
\\ \hb & \ha
\end{smallmatrix}\big)$ we can use part (1) and the result
follows.

\end{enumerate}
\end{proof}

\begin{theorem}\label{infinitesimal-induction-theorem} Let $\la=(\lx,\ly) \in \Lambda_2$ and $\mu=(\lx,m_2) \in \PP_{\la}$. Then for any $\rho \in
\hat{C}_{\mu}$ the infinitesimal inductions $\ie{\mu}{\la}(\rho)$
and $\iq{\mu}{\la}(\rho)$ are irreducible of dimension
\[
i_{\la}=\left\{
          \begin{array}{ll}
            q^{\ly-1}(q-1) , & \hbox{if $\lx > \ly$;} \\
            q^{\kk-2}(q^2-1), & \hbox{if $\lx=\ly=\kk$.}
          \end{array}
        \right.
\]
Nonequivalent representations induce to nonequivalent
representations.
\end{theorem}

\begin{proof} We argue for $\ie{\mu}{\la}(\rho)$, the assertion for $\iq{\mu}{\la}(\rho)$ is similar.
 We claim that
\begin{enumerate}

\item[(a)] $\dim \ie{\mu}{\la}(\rho)= i_{\la}$.

\item[(b)] $\ie{\mu}{\la}(\rho)_{|K_{\la}} \geqslant \big(\begin{smallmatrix} \ha & \ha
\\ \hb & \ha
\end{smallmatrix}\big)$.

\item [(c)] $\rho =\re{\mu}{\la} \ie{\mu}{\la}(\rho)$.
\end{enumerate}
Before proving these assertions note that they imply the theorem
since for any $\rho, \rho' \in \hat{C}_{\mu}$
\[
\langle \ie{\mu}{\la}(\rho), \ie{\mu}{\la}(\rho') \rangle=\langle
\rho, \re{\mu}{\la} \ie{\mu}{\la}(\rho') \rangle = \langle \rho ,
\rho' \rangle = \left\{
                  \begin{array}{ll}
                    1, & \hbox{if $\rho \simeq \rho'$;} \\
                    0, & \hbox{otherwise.}
                  \end{array}
                \right.
\]
To prove (a), we use the fact that $\dim(\rho)=q^{m_2-1}(q-1)$ (by
Theorem \ref{cuspidaltheorem}) together with
\[
[G_{\la}:\Pp_{\mu \hookrightarrow \la}]=\left\{
          \begin{array}{ll}
            q^{\ly-m_2} , & \hbox{if $\lx > \ly$;} \\
            q^{\kk-m_2-1}(q+1), & \hbox{if $\lx=\ly=\kk$}.
          \end{array}
        \right.
\]
As for (b), the representation $\ie{\mu}{\la}(\rho)$ must contain an
irreducible component which lies above the orbit of
$\big(\begin{smallmatrix} \ha & \ha
\\ \hb & \ha
\end{smallmatrix}\big) \in \hat{K}_{\la}$. This follows by recalling \eqref{ueinduction}
and the explicit form \eqref{first}, which implies that
$\big(\begin{smallmatrix} \ha & \ha
\\ \hb & \ha
\end{smallmatrix}\big) \leqslant \varphi^*\rho_{|K_{\la}}$, and since
$\varphi^*\rho \leqslant \ie{\mu}{\la}(\rho)_{|\Pp_{\mu \hookrightarrow \la}}$, we have
$\big(\begin{smallmatrix} \ha & \ha
\\ \hb & \ha
\end{smallmatrix}\big) \leqslant \ie{\mu}{\la}(\rho)_{|K_{\la}}$. Using
Proposition \ref{dimensionnilpotent} and (a), we conclude that
$\ie{\mu}{\la}(\rho)$ is irreducible.

To prove (c), observe that $\rho \leqslant \re{\mu}{\la}
\ie{\mu}{\la}(\rho)$ by adjointness of $\rb_{\hookrightarrow}$ and
$\ib_{\hookrightarrow}$. The assertion would then follow once we
show that $\dim \rho = \dim \re{\mu}{\la} \ie{\mu}{\la}(\rho)$. To
show that we apply the functors
\[
\re{(\lx,\ly-1)}{(\lx,\ly)}, \re{(\lx,\ly-2)}{(\lx,\ly-1)}, \ldots,
\re{(\lx,m_2)}{(\lx,m_2+1)}
\]
successively to $\ie{\mu}{\la}(\rho)$. Observe that for any $\mu
=(\lx,m) \le \mu'=(\lx,m') \in \PP_{\la}$, and any submodule
$M_{\mu} \subset M_{\la}$ of type $\mu$, there is a unique module
$M_{\mu'} \supset M_{\mu}$ of type $\mu'$, namely
$M_{\mu'}=M_{\mu}+M_{\la}[\pi^{m'}]$. We can therefore use Claim
\ref{associativity-infinitesimal} and write
\[
\re{(\lx,m_2)}{(\lx,m_2+1)}\circ \cdots
\circ\re{(\lx,\ly-1)}{(\lx,\ly)} = \re{\mu}{\la}.
\]
As we know by the proof of Proposition \ref{dimensionnilpotent} how the dimension
drops after each successive application of the functors
$\rb_{\hookrightarrow}$, the assertion follows.

\end{proof}


\subsection{Geometric induction in rank two}

Let $\theta=(\theta_1,\theta_2)$ be a character of $G_{(\lx)} \times
G_{(\ly)} \simeq \OO_{\lx}^{\times} \times \OO_{\ly}^{\times}$. In
this subsection we shall study the representations
\[
\begin{split}
\xi_{\theta}&=\is{(\lx),(\ly)}{(\lx,\ly)}(\theta)=\Ind_{\Pp_{(\lx),(\lx,\ly)}}^{\G_{(\lx,\ly)}}(\iota^*\theta) \\
\xi_{\theta}^{\vee}&=\is{(\ly),(\lx)}{(\lx,\ly)}(\theta)=\Ind_{\Pp_{(\ly),(\lx,\ly)}}^{\G_{(\lx,\ly)}}(\iota^*\theta)
\end{split}
\]
(recall the notations \eqref{third}, \eqref{forth} and
\eqref{geominduction}). The representations $\xi_{\theta}$ and
$\xi_{\theta}^{\vee}$ are of dimension
\begin{equation}\label{dimensiongeometric}
[G_{(\lx,\ly)}:\Pp_{(\lx),(\lx,\ly)}]=[G_{(\lx,\ly)}:\Pp_{(\ly),(\lx,\ly)}]=\left\{
                                                        \begin{array}{ll}
                                                          q^{\ly}, & \hbox{if $\lx > \ly$;} \\
                                                          q^{\kk-1}(q+1), & \hbox{if $\lx = \ly =\kk$.}
                                                        \end{array}
                                                      \right.
\end{equation}
As it happens for the groups $\G_{1^n}$, geometric induction might
be reducible or irreducible. As it turns out, the characters which
induce irreducibly are
\begin{equation}\label{Ckk}
\hat{C}_{(\lx),(\ly)}= \big\{ \theta=(\theta_i,\theta_{ii})  ~|~
\theta(u,u^{-1})=1,~ \forall u \in 1+\m^{\lx-1}\big\} \subset
\hat{G}_{(\lx)} \times \hat{G}_{(\ly)}
\end{equation}
For $\lx> \ly$ this means that $\theta_i$ does not factor through
$\OO_{\lx-1}^{\times}$ (as it is nontrivial on
$1+\m^{\lx-1}_{\lx}$), and if $\lx=\ly=\kk$ it means that
$\theta_{{i|}_{1+\m^{\kk-1}}} \ne \theta_{{ii|}_{1+\m^{\kk-1}}}$.
Let $p_{\kk,m}:\OO_{\kk} \to \OO_{m}$ denote the reduction map. The
following claim connects the infinitesimal and geometric inductions.

\begin{claim}\label{mixed} For $1 \le m < \ly$ the following equalities hold
\[
\begin{split}
&\is{(\lx),(\ly)}{(\lx,\ly)} \circ (\mathrm{id} \boxtimes
p_{\ly,m}^*) =\ie{(\lx,m)}{(\lx,\ly)} \circ \is{(\lx),(m)}{(\lx,m)}
\\&\is{(\ly),(\lx)}{(\lx,\ly)} \circ
(\mathrm{id} \boxtimes p_{\ly,m}^*)=\iq{(\lx,m)}{(\lx,\ly)} \circ
\is{(m),(\lx)}{(\lx,m)}.
\end{split}
\]

\end{claim}
\begin{proof} We argue for the first equality, the second follows
along the same lines. Unraveling definitions (e.g.
\eqref{geominduction}, \eqref{first} and \eqref{third}) we have that
\[
\begin{split}
&\is{(\lx),(\ly)}{(\lx,\ly)} \circ (\mathrm{id} \boxtimes
p_{\ly,m}^*)=\Ind_{P_2}^{G}\circ
\Inf_{P_2/U_2}^{P_2}, \quad \text{and} \\
&\ie{(\lx,m)}{(\lx,\ly)} \circ
\is{(\lx),(m)}{(\lx,m)}=\Ind_{P_1}^{G}
 \circ
\Inf_{P_1/U_1}^{P_1} \circ \Ind_{P_2/U_1}^{P_1/U_1} \circ
\Inf_{P_2/U_2}^{P_2/U_1},
\end{split}
\]
where $P_2=\Pp_{(\lx),(\lx,\ly)} \subset P_1=\Pp_{(\lx,m),(\lx,\ly)}
\subset G=\G_{(\lx,\ly)}$, and
\[
U_1=\left[\begin{matrix} 1 & \m^{\lx-\ly+m}_{\lx}
\\  & 1+\m^{m}_{\ly}
\end{matrix}\right] \subset U_2=\left[\begin{matrix} 1 & \m^{\lx-\ly}_{\lx}
\\  & 1+\m^{m}_{\ly}
\end{matrix}\right] \subset P_2.
\]
Since $U_j$ is normal in $P_j$, we can use Lemma \ref{inf-ind} and
get the desired equality.

\end{proof}

\subsubsection{Irreducible geometrically induced representations}
The following proposition is analogous to Proposition
\ref{dimensionnilpotent}.

\begin{prop}\label{dimensionsplit}Let $\la=(\lx,\ly) \in \Lambda_2$ and let $\rho$ be an irreducible
representation of $\G_{\la}$.
\begin{enumerate}

\item If $\lx>\ly >1$ and $\rho_{|K_{\la}} \geqslant \big(\begin{smallmatrix} \hat{u} & \ha
\\ \ha & \hat{z}
\end{smallmatrix}\big)$ with $\hat{u}\in \Oas$, then $\dim(\rho)=q^{\ly}$

\item If $\lx=\ly =\kk >1$ and $\rho_{|K_{\kk^2}} \geqslant \big(\begin{smallmatrix} \hat{u} & \ha
\\ \ha & \hat{z}
\end{smallmatrix}\big)$ with $\hat{u} \ne \hat{z}$, then $\dim(\rho)=q^{\kk-1}(q+1)$.
\end{enumerate}

\end{prop}

\begin{proof} \qquad

\begin{enumerate}
\item We argue by induction on $\ly$. Let $\rho$ be an irreducible
representation of $\G_{\la}$ which lies above the orbit of
$\big(\begin{smallmatrix} \hat{u} & \ha
\\ \ha & \hat{z}
\end{smallmatrix}\big) \in \hat{K}_{(\lx,\ly)}$. By twisting with a one dimensional character we may assume that $\hat{z}=0$.
Using the action \eqref{actionGkjonK} we verify that this orbit
contains the $q^2$ characters $\left\{\left(\begin{smallmatrix}
\hat{u} & c\hat{u}
\\ b\hat{u} & bc\hat{u}
\end{smallmatrix}\right)~|~ b,c \in \Oa \right\}$,
out of which $q$ survive the application of the functor
$\re{(\lx,\ly-1)}{(\lx,\ly)}$, namely those having $c=0$. It follows
that if $\xi= \re{(\lx,\ly-1)}{(\lx,\ly)}(\rho)$ then $\dim
(\xi)={\dim(\rho)}/{q}$. Using the explicit form \eqref{first} we
observe that $\langle \xi_{|V_1}, \ha_{V_1} \rangle = 0$, which
means that $\xi_{|K_{(\lx,\ly-1)}} \geqslant
\big(\begin{smallmatrix} \hat{u} & \ha
\\ \ha & \hat{z}'
\end{smallmatrix}\big)$ with the same $\hat{u} \ne 0$ and some $\hat{z}' \in \Oa$. By the induction hypothesis it follows that $\xi$
is irreducible and that $\dim(\xi)=q^{\ly-1}$, hence the assertion
follows. The basis of the induction is $\ly=2$ in which case the
role of $K$ is taken over by $H$ of \S\ref{sec:Gk1} and $\xi$ lies
above the orbits $q$-dimensional representations of the Heisenberg
group.

\item Assume $\rho \in \hat{\G}_{\kk^2}$ such that $\rho$ lies above the orbit
$\big[\big(\begin{smallmatrix} \hat{u} & \ha
\\ \ha & \hat{z}
\end{smallmatrix}\big)\big]$ with $\hat{u} \ne \hat{z}$. By twisting
with a one dimensional character we may assume that $\hat{z}=0$ and
that $\hat{u} \ne 0$. This orbit contains the $(q^2-1)$ characters
\[
\left(\begin{matrix} \frac{ad}{ad-bc}\hat{u} &
\frac{ac}{ad-bc}\hat{u}
\\ \frac{bd}{ad-bc}\hat{u} & \frac{-bc}{ad-bc}\hat{u}
\end{matrix}\right), \qquad \left(\begin{matrix} a & b \\ c & d \end{matrix}\right) \in
\GL_2(\Oa)
\]
(using the action \ref{actionGkkonK}), out of which $(q-1)$ survive
the application of the functor $\re{(\kk,\kk-1)}{(\kk,\kk)}$, namely
those having $c=0$. It follows that if $\xi=
\re{(\kk,\kk-1)}{(\kk,\kk)}(\rho)$ then $\dim
(\xi)={\dim(\rho)}/(q+1)$. Observing that $\xi_{|K_{(\kk,\kk-1)}}
\geqslant \big(\begin{smallmatrix} \hat{u} & \ha
\\ \ha & \hat{z}
\end{smallmatrix}\big)$ with $\hat{u} \ne 0$ we can use part (1) and the result
follows.

\end{enumerate}

\end{proof}

\begin{theorem} For any $\theta \in \hat{C}_{(\lx),(\ly)}$ the
geometric inductions $\xi_{\theta}$ and $\xi_{\theta}^{\vee}$ are
equivalent and irreducible. If $\lx
> \ly$ then $\theta$ determines $\xi_{\theta}$ completely. If $\lx= \ly$ then
$\xi_{\theta} \simeq \xi_{\theta'}$ if and only if $\theta=\theta'$
or $\theta^{\text{op}}=\theta'$.
\end{theorem}

\begin{proof} Since $\theta_{|V_1V_2} \leqslant \xi_{\theta|V_1V_2}$ and $\theta_{|V_1V_2}  \leqslant \xi_{\theta|V_1V_2}^{\vee}$, the
representations $\xi_{\theta}$ and $\xi_{\theta}^{\vee}$ contain the
orbits of Proposition \ref{dimensionsplit} and by dimension counting
must be irreducible. To prove the equivalence of  $\xi_{\theta}$ and
$\xi_{\theta}^{\vee}$ we shall use the identity
\[
\xi_{\theta}=\is{(\lx),(\ly)}{(\lx,\ly)}(\theta_i,\theta_{ii})={\det}^*\theta_2
\otimes \is{(\lx),(\ly)}{(\lx,\ly)}(\theta_{ii}^{-1}\theta_i,1)
\]
(and analogously for $\xi_{\theta}^{\vee}$), where we abused
notation and wrote $\theta_{ii}^{-1}$ instead of
$p_{\lx,\ly}^*(\theta_{ii}^{-1})$. Note that by \eqref{Ckk}, the
character $\theta_{ii}^{-1}\theta_i$ is nontrivial on
$1+\m^{\lx-1}_{\lx}$, hence we are reduced to the case
$\theta_{ii}=1$. In such case we also have
\[
\re{(\lx,1)}{(\lx,\ly)}(\xi_{\theta})=\rqq{(\lx,1)}{(\lx,\ly)}(\xi_{\theta})=
\re{(\lx,1)}{(\lx,\ly)}(\xi_{\theta}^{\vee})=\rqq{(\lx,1)}{(\lx,\ly)}(\xi_{\theta}^{\vee}),
\]
since all these representations are $q$-dimensional irreducible
representations of $\G_{(\lx,1)}$ with the same central character.
In particular, using Claim \ref{mixed} we get
\[
\langle \xi_{\theta}, \xi_{\theta}^{\vee} \rangle = \langle
\ie{(\lx,1)}{(\lx,\ly)}\is{(\lx),(1)}{(\lx,1)}(\theta),\xi_{\theta}^{\vee}
\rangle=\langle \is{(\lx),(1)}{(\lx,1)}(\theta),
\re{(\lx,1)}{(\lx,\ly)}(\xi_{\theta}^{\vee}) \rangle=1.
\]
Finally, for $\theta_1, \theta_2  \in \hat{C}_{(\lx),(\ly)}$ let
$\big(\begin{smallmatrix} \hat{u}_j & \ha
\\ \ha & \hat{z}_j
\end{smallmatrix}\big)$ denote the representative of
$\xi_{\theta_j|K_{\la}}$ ($j=1,2$). We should prove that
\[
\xi_{\theta_1} \simeq  \xi_{\theta_2} ~\Longleftrightarrow~ \left\{
                                                   \begin{array}{ll}
                                                     \theta_1=\theta_2, & \hbox{if $\la$ is non-rectangular;} \\
                                                     \theta_1=\theta_2 ~\mathrm{or}~ \theta_1=\theta_2^{\text{op}}
, & \hbox{if $\la$ is rectangular.}
                                                   \end{array}
                                                 \right.
\]
We argue by induction on $\ly$. If $\la$ is non-rectangular, then
$\langle \xi_{\theta_1}, \xi_{\theta_2} \rangle=0$ unless
$(\hat{u}_1,\hat{z}_1)=(\hat{u}_2,\hat{z}_2)$. If the latter holds
then we can write $\xi_{\theta_j}={\det}^*\chi_{\hat{z}} \otimes
\xi_{\theta_j'}$ with
$\theta_j'=(\chi_{\hat{z}}^{-1},\chi_{\hat{z}}^{-1})\cdot
\theta_j=(\mathrm{id} \boxtimes p_{\ly,\ly-1}^*)(\theta_j'')$ for
some $\theta_j'' \in \hat{C}_{(\lx),(\ly-1)}$. Using Claim
\ref{mixed} we now have
$\xi_{\theta_j'}=\ie{(\lx,\ly-1)}{(\lx,\ly)}(\xi_{\theta_j''})$, and
by the induction hypothesis $\langle \xi_{\theta_1''},
\xi_{\theta_2''} \rangle = \langle \theta_1'',\theta_2'' \rangle$,
which in turn implies $\theta_1=\theta_2$. If $\la$ is rectangular
then $\langle \xi_{\theta_1},\xi_{\theta_2} \rangle =0$ unless
either $(\hat{u}_1,\hat{z}_1)=(\hat{u}_2,\hat{z}_2)$ or
$(\hat{u}_1,\hat{z}_1)=(\hat{z}_2,\hat{u}_2)$. If the former holds
we proceed as in the non-rectangular case. If the latter holds,
replacing $\theta_2$ with $\theta_2^{\text{op}}$ brings us back to
the first case. The replacement does not change the representation
since the representations $\xi_{\theta}$ and
$\xi_{\theta^{\text{op}}}^{\vee}$ are equal, being intertwined by
the Weyl element. The basis for the induction is $\ly=1$ which holds
by the the remarks above.
\end{proof}

\subsubsection{Reducible geometrically induced representations} The last
ingredient of the induced representations is geometric induction
which is reducible. Let $\rho$ be one of the
$2q^{\lx-2}(q-1)\DD^{q-1} \subset \hat{\G}_{(\lx,1)}$ irreducible
representations which lie above the orbits $B_{\pm}$ of
\S\ref{sec:Gk1}. Let
\[
\xi_{\rho}=\left\{
             \begin{array}{ll}
               \ie{(\lx,1)}{(\lx,\ly)}(\rho), & \hbox{if $\rho$ lies above $\text{B}_+$;} \\
               \iq{(\lx,1)}{(\lx,\ly)}(\rho), & \hbox{if $\rho$ lies above $\text{B}_-$.}
             \end{array}
           \right.
\]

\begin{prop} The representations $\xi_{\rho}$ are irreducible of dimension $\imath_{\la}$, distinct,
and contained in geometrically induced representations.
\end{prop}

\begin{proof} The same arguments given in Theorem
\ref{infinitesimal-induction-theorem} apply to $\xi_{\rho}$ and
imply that $\dim \xi_{\rho}=\imath_{\la}$ and that it is
irreducible. Applying the functor $\rs{(\lx),(1)}{(\lx,1)}$ to
$\rho$, if $\rho$ lies above $\text{B}_+$, or the functor
$\rs{(1),(\lx)}{(\lx,1)}$ if $\rho$ lies above $\text{B}_-$, gives a
nonzero representation of $\G_{(\lx)} \times \G_{(1)}$ which proves
that $\rho$ is contained in a geometrically induced representation,
say
\[
\rho \leqslant \left\{
           \begin{array}{ll}
             \is{(\lx),(1)}{(\lx,1)}(\theta), & \hbox{if $\rho$ lies above $\text{B}_+$;} \\
             \is{(1),(\lx)}{(\lx,1)}(\theta), & \hbox{if $\rho$ lies above $\text{B}_-$.}
           \end{array}
         \right.
, \qquad \theta \in \hat{G}_{(\lx)} \times \hat{G}_{(\ly)}
\smallsetminus \hat{C}_{(\lx),(\ly)}.
\]
This means that
\[
\xi_{\rho} \leqslant \left\{
                 \begin{array}{ll}
                   \ie{(\lx,1)}{(\lx,\ly)} \circ \is{(\lx),(1)}{(\lx,1)}(\theta), & \hbox{
if $\rho$ lies above $\text{B}_+$;} \\
                   \iq{(\lx,1)}{(\lx,\ly)} \circ
\is{(1),(\lx)}{(\lx,1)}(\theta), & \hbox{if $\rho$ lies above
$\text{B}_-$,}
                 \end{array}
               \right.
\]
and the proposition now follows from Claim \ref{mixed}.
\end{proof}

\pagebreak


\section{Classification of irreducible representations}\label{sec:unitarydual}
We are now in a position to integrate all the results obtained so
far and give a complete classification of the irreducible
representations of $\G_{\la}$ $(\la \in \Lambda_2$).

\subsection{Non-rectangular case}

Let $\la=(\lx,\ly)$ with $\lx>\ly>1$. The primitive irreducibles
representations of $\G_{\la}$ which were obtained in
\S\ref{sec:constructioncuspidal} and \S\ref{sec:inductions} are
listed in the following table.

\bigskip

\centerline{\bf \underline{Table 3: Primitive representations of
$\G_{\la}$}}
\begin{align*}
& \quad \text{\em Type} &      &  \text{\em Dimension} &              & \qquad \text{\em Number}   & &\\
& \text{Cuspidal} &     & q^{\ly-1}(q-1)    &             & q^{\lx+\ly-3}(q-1)^2& &\\
&  \ie{\mu}{\la}(\hat{C}_{\mu}), ~\mu \in \{(\lx,m)\}_{m=1}^{\ly-1}
&
   &q^{\ly-1}(q-1) &   &  q^{\lx+m-3}(q-1)^2
&   & \\
&  \iq{\mu}{\la}(\hat{C}_{\mu}), ~\mu \in \{(\lx,m)\}_{m=1}^{\ly-1}
&
   &q^{\ly-1}(q-1) &   &  q^{\lx+m-3}(q-1)^2
&   & \\
& \is{(\lx),(\ly)}{(\lx,\ly)}(\hat{C}_{(\lx),(\ly)})& &\quad
q^{\ly}& & q^{\lx+\ly-3}(q-1)^3& & \\
& \xi_{\rho} < \is{(\lx),(\ly)}{(\lx,\ly)}(\hat{G}_{(\lx)} \times
\hat{G}_{(\ly)} \smallsetminus \hat{C}_{(\lx),(\ly)})& &
q^{\ly-1}(q-1) & & \quad q^{\lx-2}(q-1)& &\\& \xi_{\rho} <
\is{(\ly),(\lx)}{(\lx,\ly)}(\hat{G}_{(\lx)} \times \hat{G}_{(\ly)}
\smallsetminus \hat{C}_{(\lx),(\ly)})& & q^{\ly-1}(q-1) & & \quad
q^{\lx-2}(q-1)& &
\end{align*}

The following theorem ties up all the loose ends we have left
concerning the groups $\G_{\la}$ for $\la$ non-rectangular and
complete the classification of the irreducible representations of
$\G_{\la}$.
\begin{theorem} The representations in Table 3 exhaust all the primitive irreducible representations
of $\G_{\la}$ ($\la \in \Lambda_2$ non rectangular). In particular,
Theorem \ref{poly} part (2) and Theorem \ref{maintheorem} for $\la$
non-rectangular hold.
\end{theorem}

\begin{proof} The families to which the representations are divided
are disjoint, and the result follow by checking that we have
accumulated enough representations:
\[
|\G_{\la}|-q|\G_{\f{\la}}|=\sum \text{Number}\cdot
(\text{Dimension})^2=q^{\lx+3\ly-5}(q^3-1)(q-1)^2,
\]
where the sum runs over all the representations in Table 3.
\end{proof}

\subsection{Rectangular case} In the case of the group $\G_{\kk^2}$, the classification is different in two conceptual
points and of course in the numerics. First, the infinitesimal
functors $\ie{\mu}{\kk^2}$ and $\iq{\mu}{\kk^2}$ coincide. Second,
an infinitesimally induced representation could be further twisted
by primitive characters of $\G_{\kk^2}$ since only the orbit
$\big(\begin{smallmatrix} \ha & \hb \\ \ha & \ha
\end{smallmatrix}\big) \in \hat{K}_{\kk^2}$ can be obtained.

\bigskip

\centerline{\bf \underline{Table 4: Primitive representations of
$\G_{\kk^2}$}}
\begin{align*}
& \quad \text{\em Type} &      &  \text{\em Dimension} &              & \qquad \text{\em Number}   & &\\
& \text{Cuspidal} &     & q^{\kk-1}(q-1)    &             & \gfrac{1}{2}(q-1)(q^2-1)q^{2\kk-3}& &\\
&  q \cdot \ib_{\mu}^{\kk^2}(\hat{C}_{\mu}), \mu \in
\{(\kk,m)\}_{m=1}^{\kk-1} &
   &q^{\kk-2}(q^2-1) &   & \quad q^{\kk+m-2}(q-1)^2
&   & \\
& \is{(\kk),(\kk)}{(\kk,\kk)}(\hat{C}_{(\kk),(\kk)})& &
q^{\kk-1}(q+1)& &\quad \gfrac{1}{2}q^{2\kk-3}(q-1)^3& & \\
& \rho < \is{(\kk),(\kk)}{(\kk,\kk)}(\hat{G}_{(\kk)} \times
\hat{G}_{(\kk)} \smallsetminus \hat{C}_{(\kk),(\kk)} )& &
q^{\kk-2}(q^2-1)& & \quad q^{\kk-1}(q-1)& &
\end{align*}

\begin{theorem} The representations in Table 4 exhaust all the primitive irreducible representations
of $\G_{\kk^2}$. In particular, Theorem \ref{poly} part (3) and
Theorem \ref{maintheorem} for $\la$ rectangular hold.
\end{theorem}

\begin{proof} The families to which the representations are divided
are disjoint, and the result follow by checking that we have
accumulated enough representations:
\[
|\G_{\kk^2}|-q|\G_{\f{\kk^2}}|=\sum \text{Number}\cdot
(\text{Dimension})^2=q^{4\kk-7}(q^3-1)(q^2-q)(q^2-1),
\]
where the sum runs over all the representations in Table 4.
\end{proof}

To make a link with the literature we remark that the
representations of type (ii) were constructed in \cite{Nobs1} where
they are called {\em la s\myacute{e}rie
d\myacute{e}ploy\myacute{e}e}. The cuspidal representations were
constructed there as well and are called {\em la s\myacute{e}rie non
ramifi\myacute{e}e}. The infinitesimally induced representations
turn out to be the missing part in the description of the
irreducible representations of $\G_{\kk^2}$, and was the prime
motivation behind the present study. Note that the
infinitesimal induction from cuspidal representations of $\G_{\mu}$
($\mu \in \PP_{\kk^2}$) should be further twisted in order to cover
all the orbits, as only the orbit with $\hat{u}=0$ in $(iii)$ of
Table 2 survives after an application of infinitesimal restriction
functors.

It is perhaps worthwhile mentioning that though a complete and
unified description of the rectangular case seems to be missing in
the literature, for this particular case all the necessary
ingredients are handy and one could complete the
classification, see a recent preprint by Stasinski \cite{stasinski-2008}.

Another remark which is in order regards even characteristic. In
\cite{Nobs1.5,Nobs2} a separate discussion is devoted to the even
characteristic case, and even that only for $\GL_2(\Z_2)$. It is
important to note that the difference lies in the method and not in
the objects under study. This has been just shown to hold in
general: for any $\G_{\la}$ ($\la \in \Lambda_2$), and over any
$\OO$, the cardinality of the residue field $\OO/\m$ appears only as a
parameter and not in any essential way.


\section{Conjugacy classes}\label{sec:conjclasses}

Recall that for any finite $\OO$-module $M$, an element $f \in
\End_{\OO}(M)$ is called cyclic if there exist $m \in M$ such that
$m, f \cdot m, f^2 \cdot m, \ldots ,f^j \cdot m$ span $M$ over $\OO$
for some $j \in \N$.

\begin{defn} An element $f \in \End_{\OO}(M)$ is called {\em almost cyclic} if
it can be written as
\[
f=a I + \pi^{i}h
\]
with $a \in \OO$, $i$ a nonnegative integer and $h \in
\End_{\OO}(M)$ a cyclic element
 with $\pi^ih \ne 0$.
\end{defn}
The classification of conjugacy classes of $\G_{\la}$ ($\la \in
\Lambda_2$) is greatly simplified by the following dichotomy:

\begin{lemma} Any element in $E_{\la}=\End_{\OO}(M_{\la})$ ($\la \in
\Lambda_2$) is either scalar or almost cyclic.
\end{lemma}

\begin{proof}
Let $f \in E_{\la}$ be nonscalar. Then $f=aI+\pi^i h$ for some
nonnegative integer $i$, $a \in \OO$, and $h \in E_{\la}$ such that
$\pi^ih \ne 0$ and $\bar{h} \in \End_{\OO}(M_{\la} /\pi M_{\la})$ is
nonscalar. This implies that $\bar{h}$ is cyclic, since $M_{\la}/\pi
M_{\la} \simeq \OO_1^2$, and any element in $\End_{\OO}(\OO_1^2)$ is
either scalar or cyclic. By Nakayama's lemma $h$ must be cyclic
since it is a lift of a cyclic element.

\end{proof}

We shall now build the conjugacy classes of $\G_{\la}$ inductively.
The inverse image $p^{-1}(C)$ of a conjugacy class $C \subset
\G_{\f{\la}}$ with respect to the reduction map $p:\G_{\la} \to
\G_{\f{\la}}$ is a disjoint union of conjugacy classes in
$\G_{\la}$. An almost cyclic class has precisely $q^2$ classes above
it while the classes in $\G_{\la}$ which lie above a scalar class
$aI \in \G_{\f{\la}}$ are parameterized by the $\G_{\la}$-orbits in
$p^{-1}(aI)=aI \cdot K_{\la} \subset \G_{\la}$ under conjugation.
This is given by
\[
\begin{split}
&\#\text{Orbits of $\G_{\kk^2}$ on $K_{\kk^2}$}:
\quad q^2+q, \qquad \kk > 1 \\
&\#\text{Orbits of $\G_{\la}$ on $K_{\la}$}:
\quad q^2+q+1, \qquad \lx > \ly > 1 \\
\end{split}
\]
where the first equality is simply the number of similarity classes
in $K_{\kk^2} \simeq M_2(\OO_1)$: $q$ scalar and $q^2$ cyclic. The
second equality follows by recalling the action \eqref{action} in
the special case $K_{\la}=K_{\la}^{1,0} \simeq M_2(\Oa)$. Explicitly
the orbits are represented in this case by
\[
\left\{\left(\begin{matrix} u & 0 \\ 0 & 0
\end{matrix}\right) ~\Big|~ u \in \Oa\right\}, \left\{\left(\begin{matrix} 0 & v \\ 1 & z
\end{matrix}\right)~\Big|~ v,z \in \Oa \right\}, \left(\begin{matrix} 0 & 1 \\ 0 & 0
\end{matrix}\right).
\]
With this in hand we can prove the following.

\begin{theorem}
For any $\la \in \Lambda_2$ the number of conjugacy classes in
$\G_{\la}$ is
\begin{equation}\label{numberconj}
|\Conj(\G_{\la})| = \left\{
                      \begin{array}{ll}
q^{2\kk}-q^{\kk-1}, & \hbox{if $\lx = \ly =\kk$;}
\\
 q^{\lx+\ly-2}(q^2-q+2)-q^{\lx-2}(q+1), & \hbox{if $\lx > \ly$.}
                      \end{array}
                    \right.
\end{equation}
\end{theorem}

\begin{proof} We argue by induction on the type. The basis for the
induction is given by $\G_{(\kk,1)}$ ($\kk \ge 1$):
\[
\begin{split}
&\Conj(\G_{(1,1)}) \simeq \left\{\left(\begin{smallmatrix} u & 0 \\
0 & u
\end{smallmatrix}\right) | u \in \Oas\right\} \sqcup \left\{\left(\begin{smallmatrix} 0 & v \\ 1 & z
\end{smallmatrix}\right)  | v \in \Oa,~z \in \Oas \right\}, \\
&\Conj(\G_{(\kk,1)}) \simeq \left\{\left(\begin{smallmatrix} u & 0
\\ 0 & z
\end{smallmatrix}\right) | u,z \in \Oks \right\} \sqcup \left\{\left(\begin{smallmatrix} u &
\pi^{\kk-1}
\\ 0 & u
\end{smallmatrix}\right) | u \in \OO_{\kk-1}^{\times} \right\} \sqcup \left\{\left(\begin{smallmatrix} u &
w
\\ 1 & u
\end{smallmatrix}\right) | u \in \OO_{\kk-1}^{\times},~w \in \Oa \right\}, \quad
\text{if $\kk >1 $}.
\end{split}
\]
 In particular, $|\Conj(\G_{(1,1)})|=q^2-1$ and
$|\Conj(\G_{(\kk,1)})|=q^3-q^2+q-1$ for $\kk>1$, which establishes
the basis for the induction. Assuming \eqref{numberconj} for $\la$
we shall prove it for $\sharp{\la}=(\lx+1,\ly+1)$. Indeed,
$\Conj(\G_{\la})$ contains $|\Oks|$ scalar elements and the rest are
almost cyclic. Lifting them to $\G_{\sharp{\la}}$ gives
\[
\begin{split}
|\Conj(\G_{\sharp{\la}})|&=(q^2+q)|\Oks|+q^2(|\Conj(\G_{\la})|-|\Oks|)\\
&=(q^2+q)q^{\kk-1}(q-1)+q^2[(q^{2\kk}-q^{\kk-1})-q^{\kk-1}(q-1)]\\
&=q^{2\kk+2}-q^{\kk}
\end{split}
\]
if $\la=\kk^2$, and
\[
\begin{split}
|\Conj(\G_{\sharp{\la}})|&=(q^2+q+1)|\OO_{\lx}^{\times}|+q^2(|\Conj(\G_{\la})|-|\OO_{\lx}^{\times}|)\\
&=q^{\lx+\ly}(q^2-q+2)-q^{\lx-1}(q+1)
\end{split}
\]
if $\lx > \ly >1$.

\end{proof}

The result of course matches the polynomials defined by the
representations (see Theorem \ref{poly}), that is
$R_{\la}(1)=|\Conj(\G_{\la})|$.

\vspace{\bigskipamount}

\begin{footnotesize}
\begin{quote}

Uri Onn\\
Department of Mathematics, Ben-Gurion university of the Negev,\\
Beer-Sheva 84105, Israel\\
{\tt urionn@math.bgu.ac.il}

\end{quote}
\end{footnotesize}


\begin{thebibliography}{10}

\bibitem{AOPS}
{\sc A.-M. Aubert, U.~Onn, A.~Prasad, and A.~Stasinski}, {\em On cuspidal
  representations of general linear groups over discrete valuation rings}, to
  appear in the Israel Journal of Mathematics, arXiv:0706.0058 [math.RT].

\bibitem{BO1}
{\sc U.~Bader and U.~Onn}, {\em {G}eometric representations of
  $\mathrm{GL}(n,\mathrm{R})$, cellular {H}ecke algebras and the embedding
  problem}, Journal of Pure and Applied Algebra, 208 (2007), pp.~905--922.

\bibitem{carayol}
{\sc H.~Carayol}, {\em Repr\'esentations cuspidales du groupe lin\'eaire}, Ann.
  Sci. \'Ecole Norm. Sup. (4), 17 (1984), pp.~191--225.

\bibitem{Casselman}
{\sc W.~Casselman}, {\em The restriction of a representation of {${\rm
  GL}\sb{2}(k)$} to {${\rm GL}\sb{2}({\mathcal O})$}}, Math. Ann., 206 (1973),
  pp.~311--318.

\bibitem{Gerardin}
{\sc P.~Gerardin}, {\em {Construction de s\'eries discretes p-adiques. Sur les
  s\'eries discretes non ramifiees des groupes reductifs deployes p-adiques.}},
  {Lecture Notes in Mathematics. 462. Springer-Verlag}, 1975.

\bibitem{Green}
{\sc J.~A. Green}, {\em The characters of the finite general linear groups},
  Trans. Amer. Math. Soc., 80 (1955), pp.~402--447.

\bibitem{HC}
{\sc Harish-Chandra}, {\em Eisenstein series over finite fields}, in Functional
  analysis and related fields (Proc. Conf. M. Stone, Univ. Chicago, Chicago,
  Ill., 1968), Springer, New York, 1970, pp.~76--88.

\bibitem{HG1}
{\sc G.~Hill}, {\em A {J}ordan decomposition of representations for {${\rm
  GL}_n(\mathcal O)$}}, Comm. Algebra, 21 (1993), pp.~3529--3543.

\bibitem{HG2}
\leavevmode\vrule height 2pt depth -1.6pt width 23pt, {\em On the nilpotent
  representations of {${\rm GL}_n(\mathcal O)$}}, Manuscripta Math., 82 (1994),
  pp.~293--311.

\bibitem{HG4}
\leavevmode\vrule height 2pt depth -1.6pt width 23pt, {\em Regular elements and
  regular characters of {${\rm GL}_n({\mathcal O})$}}, J. Algebra, 174 (1995),
  pp.~610--635.

\bibitem{HG3}
\leavevmode\vrule height 2pt depth -1.6pt width 23pt, {\em Semisimple and
  cuspidal characters of {${\rm GL}_n({\mathcal O})$}}, Comm. Algebra, 23
  (1995), pp.~7--25.

\bibitem{Isaacs}
{\sc I.~M. Isaacs}, {\em Character theory of finite groups}, Academic Press
  [Harcourt Brace Jovanovich Publishers], New York, 1976.
\newblock Pure and Applied Mathematics, No. 69.

\bibitem{JZ}
{\sc A.~Jaikin-Zapirain}, {\em Zeta function of representations of compact
  {$p$}-adic analytic groups}, J. Amer. Math. Soc., 19 (2006), pp.~91--118.

\bibitem{Kloosterman2}
{\sc H.~D. Kloosterman}, {\em The behaviour of general theta functions under
  the modular group and the characters of binary modular congruence groups.
  {I}}, Ann. of Math. (2), 47 (1946), pp.~317--375.

\bibitem{Kloosterman1}
\leavevmode\vrule height 2pt depth -1.6pt width 23pt, {\em The behaviour of
  general theta functions under the modular group and the characters of binary
  modular congruence groups. {II}}, Ann. of Math. (2), 47 (1946), pp.~376--447.

\bibitem{Kutzko}
{\sc P.~C. Kutzko}, {\em The characters of the binary modular congruence
  group}, Bull. Amer. Math. Soc., 79 (1973), pp.~702--704.

\bibitem{Kutzko1}
{\sc P.~C. Kutzko}, {\em On the supercuspidal representations of {${\rm
  Gl}\sb{2}$}}, Amer. J. Math., 100 (1978), pp.~43--60.

\bibitem{Lusztig1}
{\sc G.~Lusztig}, {\em Representations of reductive groups over finite rings},
  Represent. Theory, 8 (2004), pp.~1--14.

\bibitem{MI1}
{\sc I.~G. Macdonald}, {\em Symmetric functions and {H}all polynomials}, Oxford
  Mathematical Monographs, The Clarendon Press Oxford University Press, New
  York, second~ed., 1995.
\newblock With contributions by A. Zelevinsky, Oxford Science Publications.

\bibitem{Nagornyi1}
{\sc S.~V. Nagorny{\u\i}}, {\em Complex representations of the group {${\rm
  GL}(2, Z/p\sp{n}Z)$}}, Zap. Nau\v cn. Sem. Leningrad. Otdel. Mat. Inst.
  Steklov. (LOMI), 64 (1976), pp.~95--103, 161.
\newblock Rings and modules.

\bibitem{Nobs3}
{\sc A.~Nobs}, {\em Die irreduziblen {D}arstellungen der {G}ruppen
  {$SL\sb{2}(Z\sb{p})$}, insbesondere {$SL\sb{2}(Z\sb{2})$}. {I}}, Comment.
  Math. Helv., 51 (1976), pp.~465--489.

\bibitem{Nobs1}
\leavevmode\vrule height 2pt depth -1.6pt width 23pt, {\em La s\'erie
  d\'eploy\'ee et la s\'erie non ramifi\'ee de {$\rm {GL}_{2}({O})$}}, C. R.
  Acad. Sci. Paris S\'er. A-B, 283 (1976), pp.~Aii, A297--A300.

\bibitem{Nobs1.5}
\leavevmode\vrule height 2pt depth -1.6pt width 23pt, {\em Les
  repr\'esentations irr\'eductibles du groupe {${\rm GL}\sb{2}({\bf
  Z}\sb{2})$}}, C. R. Acad. Sci. Paris S\'er. A-B, 283 (1976), pp.~Ai,
  A433--A436.

\bibitem{Nobs2}
\leavevmode\vrule height 2pt depth -1.6pt width 23pt, {\em Die irreduziblen
  {D}arstellungen von {${\rm GL}\sb{2}(Z\sb{p})$}, insbesondere {${\rm
  GL}\sb{2}(Z\sb{2})$}}, Math. Ann., 229 (1977), pp.~113--133.

\bibitem{NobsWolfart1}
{\sc A.~Nobs and J.~Wolfart}, {\em Darstellungen von {${\rm SL}(2,$} {${\bf
  Z}/p\sp{\lambda }{\bf Z})$} und {T}hetafunktionen. {I}}, Math. Z., 138
  (1974), pp.~239--254.

\bibitem{NobsWolfart2}
\leavevmode\vrule height 2pt depth -1.6pt width 23pt, {\em Die irreduziblen
  {D}arstellungen der {G}ruppen {$SL\sb{2}(Z\sb{p})$}, insbesondere
  {$SL\sb{2}(Z\sb{p})$}. {II}}, Comment. Math. Helv., 51 (1976), pp.~491--526.

\bibitem{PD1}
{\sc D.~Prasad}, {\em Invariant forms for representations of {${\rm GL}\sb 2$}
  over a local field}, Amer. J. Math., 114 (1992), pp.~1317--1363.

\bibitem{Shalika}
{\sc J.~A. Shalika}, {\em Representation of the two by two unimodular group
  over local fields}, in Contributions to automorphic forms, geometry, and
  number theory, Johns Hopkins Univ. Press, Baltimore, MD, 2004, pp.~1--38.

\bibitem{Silbereger2}
{\sc A.~J. Silberger}, {\em {${\rm PGL}\sb{2}$} over the {$p$}-adics: its
  representations, spherical functions, and {F}ourier analysis}, Lecture Notes
  in Mathematics, Vol. 166, Springer-Verlag, Berlin, 1970.

\bibitem{Springer}
{\sc T.~A. Springer}, {\em Kloosterman's work on representations of finite
  modular groups}, Nieuw Arch. Wiskd. (5), 1 (2000), pp.~126--129.

\bibitem{stasinski-2008}
{\sc A.~Stasinski}, {\em The smooth representations of ${GL}_2(\mathcal{O})$},
  2008, arXiv:0807.4684v1 [math.RT].

\bibitem{Tanaka}
{\sc S.~Tanaka}, {\em Irreducible representations of the binary modular
  congruence groups {${\rm mod}\,p\sp{\lambda }$}}, J. Math. Kyoto Univ., 7
  (1967), pp.~123--132.

\bibitem{Yoshida}
{\sc T.~Yoshida}, {\em On the unit groups of {B}urnside rings}, J. Math. Soc.
  Japan, 42 (1990), pp.~31--64.

\bibitem{Zelevinsky}
{\sc A.~V. Zelevinsky}, {\em Representations of finite classical groups},
  vol.~869 of Lecture Notes in Mathematics, Springer-Verlag, Berlin, 1981.
\newblock A Hopf algebra approach.

\end{thebibliography}
\end{document}